\documentclass{amsart}
\usepackage{amssymb}
\usepackage[mathcal]{euscript}
\usepackage{fullpage}
\usepackage{epsfig}

\input xypic

\newcommand{\CB}        {{\mathcal{B}}}
\newcommand{\DD}        {{\mathcal{D}}}
\newcommand{\CL}        {{\mathcal{L}}}
\newcommand{\CR}        {{\mathcal{R}}}
\newcommand{\FF}        {{\mathbb{F}}}
\newcommand{\Fp}        {{\mathbb{F}_p}}          
\newcommand{\Zh}        {{\mathbb{Z}}}
\newcommand{\Zpl}       {{\mathbb{Z}_{(p)}}}      
\newcommand{\Ztl}       {{\mathbb{Z}_{(2)}}}      
\newcommand{\real}      {{\mathbb{R}}}              
\newcommand{\cplx}      {{\mathbb{C}}}              
              
\newcommand{\rinf}      {{\mathbb{R}}^\infty}

\newcommand{\LL}        {{\mathbb{L}}}

\newcommand{\CM}        {{\mathcal{M}}}
\newcommand{\CS}        {{\mathcal{S}}}
\newcommand{\CT}        {{\mathcal{T}}}
\newcommand{\CU}        {{\mathcal{U}}}
\newcommand{\CV}        {{\mathcal{V}}}

\newcommand{\hM}        {h{\mathcal{M}}}
\newcommand{\hT}        {h{\mathcal{T}}}
\newcommand{\hbM}       {\overline{h}{\mathcal{M}}}
\newcommand{\hbT}       {\overline{h}{\mathcal{T}}}

\newcommand{\al}        {\alpha}
\newcommand{\bt}        {\beta} 
\newcommand{\btb}       {\overline{\beta}}
\newcommand{\btt}       {\beta\Smash\beta}
\newcommand{\gm}        {\gamma}
\newcommand{\dl}        {\delta}
\newcommand{\ep}        {\epsilon}
\newcommand{\tht}       {\theta}
\newcommand{\om}        {\omega}
\newcommand{\lm}        {\lambda}
\newcommand{\sg}        {\sigma}

\newcommand{\Der}       {\operatorname{Der}}
\newcommand{\Hom}       {\operatorname{Hom}}

\newcommand{\cpi}       {{\mathbb{C}P^\infty}}
\newcommand{\rp}        {{\mathbb{R}P}} 
\newcommand{\rpi}       {{\mathbb{R}P^\infty}}
\newcommand{\tm}        {\times}
\newcommand{\Rxx}       {(R/x)^{(2)}}
\newcommand{\Rxxx}      {(R/x)^{(3)}}
\newcommand{\BP}[1]     {BP\langle #1\rangle}
\newcommand{\Sg}        {\Sigma}
\newcommand{\Smash}     {\wedge}
\newcommand{\ann}       {\operatorname{ann}}
\newcommand{\cb}        {\overline{c}}
\newcommand{\cp}        {{\mathbb{C}P}} 
\newcommand{\sse}       {\subseteq}
\newcommand{\st}        {\;|\;}
\newcommand{\tP}        {\widetilde{P}}
\newcommand{\xra}       {\xrightarrow}
\newcommand{\Om}        {\Omega}
\newcommand{\Omi}       {\Omega^\infty}
\newcommand{\Sgi}       {\Sigma^\infty}
\newcommand{\ov}[1]     {\overline{#1}}
\newcommand{\Qb}        {\overline{Q}}
\newcommand{\ths}       {\ltimes}
\newcommand{\fps}[2]    {{#1 [\![ #2 ]\!]}} 
\newcommand{\Wedge}     {\vee}
\newcommand{\half}      {{\textstyle\frac{1}{2}}}
\newcommand{\sixth}     {{\textstyle\frac{1}{6}}}
\newcommand{\xla}       {\xleftarrow}
\newcommand{\mra}       {\xrightarrow{}}
\newcommand{\era}       {\xrightarrow{}}
\newcommand{\pb}        {\overline{p}}
\newcommand{\xb}        {\overline{x}}
\newcommand{\lb}        {\left[}
\newcommand{\rb}        {\right]}
\newcommand{\ot}        {\otimes}
\newcommand{\even}      {{\text{ev}}}
\newcommand{\eqdef}     {:=}

\newcommand{\invlim} {\operatornamewithlimits{\underset{\longleftarrow}{lim}}}

\newcommand{\hcolim} {\operatornamewithlimits{\underset{\longrightarrow}{holim}}}

\renewcommand{\:}{\colon}

\newtheorem{theorem}{Theorem}[section]

\newtheorem{lemma}[theorem]{Lemma}
\newtheorem{proposition}[theorem]{Proposition}
\newtheorem{corollary}[theorem]{Corollary}
\theoremstyle{definition}
\newtheorem{remark}[theorem]{Remark}
\newtheorem{definition}[theorem]{Definition}

\newtheorem{construction}[theorem]{Construction}

\begin{document}
\title{Products on $MU$-modules}
\author{N.~P.~Strickland}
\address{Trinity College, Cambridge CB2 1TQ, England}
\email{neil@dpmms.cam.ac.uk}
\bibliographystyle{abbrv}

\begin{abstract}
 Elmendorf, Kriz, Mandell and May have used their technology of
 modules over highly structured ring spectra to give new constructions
 of $MU$-modules such as $BP$, $K(n)$ and so on, which makes it much
 easier to analyse product structures on these spectra.
 Unfortunately, their construction only works in its simplest form for
 modules over $MU[\half]_*$ that are concentrated in degrees divisible
 by $4$; this guarantees that various obstruction groups are trivial.
 We extend these results to the cases where $2=0$ or the homotopy
 groups are allowed to be nonzero in all even degrees; in this context
 the obstruction groups are nontrivial.  We shall show that there are
 never any obstructions to associativity, and that the obstructions to
 commutativity are given by a certain power operation; this was
 inspired by parallel results of Mironov in Baas-Sullivan theory.  We
 use formal group theory to derive various formulae for this power
 operation, and deduce a number of results about realising $2$-local
 $MU_*$-modules as $MU$-modules.
\end{abstract}

\maketitle 

\section{Introduction}

A great deal of work in algebraic topology has exploited the
generalised cohomology theory $MU^*(X)$ (for spaces $X$), which is
known as complex cobordism; good entry points to the literature
include~\cite{ad:shg,ra:nps,ra:ccs,wi:bph}.  This theory is
interesting because of its connection with the theory of formal group
laws (FGL's), starting with Quillen's fundamental
theorem~\cite{qu:fgl,qu:epc} that $MU^*$ is actually the universal
example of a ring equipped with an FGL.

Suppose that we have a graded ring $A^*$ equipped with an FGL.  In the
cases discussed below, the FGL involved will generally be the
universal example of an FGL with some interesting property.  Examples
include the rings known to topologists as $BP^*$, $P(n)^*$, $K(n)^*$
and $E(n)^*$; see Section~\ref{sec-statements} for the definitions.
It is natural to ask whether there is a generalised cohomology theory
$A^*(X)$ whose value on a point is the ring $A^*$, and a natural
transformation $MU^*(X)\xra{}A^*(X)$, such that the resulting map
$MU^*\xra{}A^*$ carries the universal FGL over $MU^*$ to the given FGL
over $A^*$.  This question has a long history, and has been addressed
by a number of different methods for different rings $A^*$.  The
simplest case is when $A^*$ is obtained from $MU^*$ by inverting some
set $S$ of nonzero homogeneous elements, in other words
$A^*=S^{-1}MU^*$.  In that case the functor
$A^*(X)=A^*\ot_{MU^*}MU^*(X)$ is a generalised cohomology theory on
finite complexes, which can be extended to infinite complexes or
spectra by standard methods.  For example, given a prime $p$ one can
invert all other primes to get a cohomology theory $MU_{(p)}^*(X)$.
Cartier had previously introduced the notion of a $p$-typical FGL and
constructed the universal example of such a thing over $BP^*$, which
is a polynomial algebra over $\Zpl$ on generators $v_k$ in degree
$-2(p^k-1)$ for $k>0$.  It was thus natural to ask our ``realisation
question'' for $A^*=BP^*$.  Quillen~\cite{qu:fgl} constructed an
idempotent self map $\ep\:MU_{(p)}^*(X)\xra{}MU_{(p)}^*(X)$, whose
image is a subring, which we call $BP^*(X)$.  He showed that this is a
cohomology theory whose value on a point is the ring $BP^*$, and that
the FGL's are compatible in the required manner.  This cohomology
theory was actually defined earlier by Brown and
Peterson~\cite{brpe:sca} (hence the name), but in a less structured
and precise way.  It is not hard to check that we again have
$BP^*(X)=BP^*\ot_{MU^*}MU^*(X)$ when $X$ is finite.  This might tempt
us to just define $A^*(X)=A^*\ot_{MU^*}MU^*(X)$ for any $A^*$, but
unfortunately this does not usually have the exactness properties
required of a generalised cohomology theory.  Another major advance
was Landweber's determination~\cite{la:hpc} of the precise conditions
under which $A^*\ot_{MU^*}MU^*(X)$ does have the required exactness
properties, which turned out to be natural ones from the point of view
of formal groups.  However, there are many cases of interest in which
Landweber's exactness conditions are not satisfied, and for these
different methods are required.  Many of them are of the form
$A^*=(S^{-1}MU^*)/I$ for set $S$ of homogeneous elements and some
homogeneous ideal $I\leq S^{-1}MU^*$.  For technical reasons things
are easier if we assume that $I$ is generated by a regular sequence,
in other words $I=(x_1,x_2,\ldots)$ and $x_k$ is not a zero-divisor in
$(S^{-1}MU^*)/(x_j\st j<k)$.  If $A^*$ arises in this way, we say that
it is a \emph{localised regular quotient (LRQ)} of $MU^*$.  If
$S=\emptyset$ we say that $A^*$ is a \emph{regular quotient} of
$MU^*$.  The first advance in this context was the Baas-Sullivan
theory of cobordism of manifolds with singularities~\cite{ba:bms}.
Given a regular quotient $A^*$ of $MU^*$, this theory constructed a
cohomology theory $A^*(X)$, landing in the category of $MU^*$-modules,
and a map $MU^*(X)\xra{}A^*(X)$.  Unfortunately, the details were
technically unwieldy, and it was not clear whether $A^*(X)$ was unique
or whether it had a natural product structure, and if so whether it
was commutative or associative.  Some of these questions were
addressed by Shimada and Yagita~\cite{shya:mcb}, Mironov~\cite{mi:ems}
and Morava~\cite{mo:pop}, largely using the geometry of cobordisms.
Another idea was (in special cases, modulo some technical details) to
calculate the group of all natural transformations $A^*(X)\ot
A^*(X)\xra{}A^*(X)$ and then see which of them are commutative,
associative and unital.  This was the approach of
W\"{u}rgler~\cite{wu:crs,wu:ctu,wu:pfc}; much more recently, Nassau
has corrected some inaccuracies and extended these
results~\cite{na:nks,na:spp}.

Baas-Sullivan theory eventually yielded satisfactory answers for rings
of the form $MU^*/x$, but the work involved in handling ideals with
more than one generator remained rather hard.  The picture changed
dramatically with the publication of~\cite{ekmm:rma} by Elmendorf,
Kriz, Mandell and May (hereafter referred to as EKMM), which we now
explain.  Firstly, the natural home for our investigation is not
really the category of generalised cohomology theories, but rather
Boardman's homotopy category of spectra~\cite{ad:shg,ma:ssa}, which we
call $\CB$.  There is a functor $\Sgi$ from finite complexes to $\CB$,
and any cohomology theory $A^*(X)$ on finite complexes is represented
by a spectrum $A\in\CB$ in the sense that $A^n(X)=[\Sgi X,\Sg^nA]$ for
all $n$ and $X$.  The representing spectrum $A$ is unique up to
isomorphism~\cite{br:ct,ad:vbr}, and the isomorphism is often unique.
There have been many different constructions of categories equivalent
to $\CB$.  The starting point of~\cite{ekmm:rma} was EKMM's
construction of a topological model category $\CM$ with a symmetric
monoidal smash product, whose homotopy category is equivalent to
$\CB$.  This was previously feared to be impossible, for subtle
technical reasons~\cite{le:ccs}.  EKMM were also able to construct a
version of $MU$ which was a strictly commutative monoid in $\CM$,
which allowed them to define the category $\CM_{MU}$ of $MU$-modules.
They showed how to make this into topological model category, and thus
defined an associated homotopy category $\DD_{MU}$.  This again has a
symmetric monoidal smash product, which should be thought of as a sort
of tensor product over $MU$.  They showed that the problem of
realising LRQ's of $MU^*$ becomes very much easier if we work in
$\DD_{MU}$ (and then apply a forgetful functor to $\CB$ if required).
In fact their methods work when $MU$ is replaced by any strictly
commutative monoid $R$ in $\CM$ such that $R^*$ is concentrated in
even degrees.  They show that if $A^*$ is an LRQ of $R^*$ and $2$ is
invertible in $A^*$ and $A^*$ is concentrated in degrees divisible by
$4$, then $A$ can be realised as a commutative and associative ring
object in $\DD_R$.

In the present work, we will start by sharpening this slightly.  The
main point here is that EKMM notice an obstruction to associativity in
$A^{4k+2}$, so they assume that these groups are zero.  Motivated by a
parallel result in Baas-Sullivan theory~\cite{mi:mtc}, we show that
the associativity obstructions are zero even if the groups are not
(see Remark~\ref{rem-ass}).  We deduce that if $A^*$ is an LRQ of
$R^*$ and $2$ is invertible in $A^*$ then $A$ can be realised as a
commutative and associative ring in $\DD_R$, in a way which is unique
up to unique isomorphism (Theorem~\ref{thm-odd-general}).  We also
prove a number of subsidiary results about the resulting ring objects.

The more substantial part of our work is the attempt to remove the
condition that $2$ be invertible in $A^*$, without which the results
become somewhat more technical.  We show that the obstruction to
defining a commutative product on $R/x$ is given by $\tP(x)$ for a
certain power operation $\tP\:R^d\xra{}R^{2d+2}/2$.  This was again
inspired by a parallel result of Mironov~\cite{mi:mtc}.  We deduce
that if $A^*=S^{-1}R^*/I$ is an LRQ of $R^*$ without $2$-torsion and
$\tP(I)\leq I\pmod{2}$ then $A^*$ is again uniquely realisable
(Theorem~\ref{thm-even-general}).  When $A^*$ has $2$-torsion we have
no such general result and must proceed case by case.  Again following
Mironov, we show that when $R=MU$, the operation $\tP$ can be computed
using formal group theory.  We considerably extend and sharpen
Mironov's calculations, using techniques which I hope will be useful
in more general work on power operations.  Using these results, we
show that many popular LRQ's of $MU_{(2)}^*$ have almost unique
realisations as associative, almost commutative rings in $\DD_{MU}$.
See Theorems~\ref{thm-MUIn} and~\ref{thm-Pn} for precise statements.
The major exceptions are the rings $\BP{n}^*$ and $E(n)^*$, but we
show that even these become uniquely realisable as commutative rings
in $\DD_{MU}$ if we allow ourselves to modify the usual definition
slightly.  We call the resulting spectra $\BP{n}'$ and $E(n)'$; they
are acceptable substitutes for $\BP{n}$ and $E(n)$ in almost all
situations.

\section{Statement of Results}
\label{sec-statements}

We use the category $\CM$ of $S$-modules as constructed
in~\cite{ekmm:rma}; we recall some details in Section~\ref{sec-found}.
The main point is that $\CM$ is a symmetric monoidal category with a
closed model structure whose homotopy category is Boardman's homotopy
category of spectra.  We shall refer to the objects of $\CM$ simply as
spectra.

Because $\CM$ is a symmetric monoidal category, it makes sense to talk
about strictly commutative ring spectra; these are essentially
equivalent to $E_\infty$ ring spectra in earlier foundational
settings.  Let $R$ be such an object, such that $R_*=\pi_*R$ is even
(by which we mean, concentrated in even degrees).  We also assume that
$R$ is $q$-cofibrant in the sense of~\cite[Chapter VII]{ekmm:rma} (if
not, we replace $R$ by a weakly equivalent cofibrant model).  The main
example of interest to us is $R=MU$ .  There are well-known
constructions of $MU$ as a spectrum in the earlier sense of Lewis and
May~\cite{lemast:esh}, with an action of the $E_\infty$ operad of
complex linear isometries.  Thus, the results 
of~\cite[Chapter II]{ekmm:rma} allow us to construct $MU$ as a
strictly commutative ring spectrum.

One can define a category $\CM_R$ of $R$-modules in the evident way,
with all diagrams commuting at the geometric level.  After inverting
weak equivalences, we obtain a homotopy category $\DD=\DD_R$, referred
to as the derived category of $\CM_R$.  We shall mainly work in this
derived category, and the category $\CR=\CR_R$ of ring objects in
$\DD$ (referred to in~\cite{ekmm:rma} as $R$-ring spectra).  All our
ring objects are assumed to be associative and to have a two-sided
unit.  Thus, an object $A\in\CR$ has an action $R\Smash_SA\xra{}A$
which makes various diagrams commute at the geometric level, and a
product $A\Smash_SA\xra{}A$ that is geometrically compatible with the
$R$-module structures, and is homotopically associative and unital.
We also write $\CR_*$ for the category of algebras over the discrete
ring $R_*$.  We write $\CR_*^e$ for the category of even
$R_*$-algebras, and $\CR_*^c$ for the commutative ones, and similarly
$\CR_*^{ec}$, $\CR^e$, $\CR^c$ and $\CR^{ec}$.

\begin{definition}
 Let $A_*$ be an even commutative $R_*$-algebra without $2$-torsion.
 A \emph{strong realisation} of $A_*$ is a commutative ring object
 $A\in\CR^{ec}$ with a given isomorphism $\pi_*(A)\simeq A_*$, such
 that the resulting map
 \[ \CR(A,B) \xra{} \CR(A_*,\pi_*(B)) \]
 is an isomorphism whenever $B\in\CR^{ec}$ and $B_*$ has no
 $2$-torsion.  We say that $A_*$ is \emph{strongly realisable} if such
 a realisation exists.
\end{definition}

\begin{remark}
 It is easy to see that the category of strongly realisable
 $R_*$-algebras is equivalent to the category of those $A\in\CR^{ec}$
 for which $\pi_*(A)$ is strongly realisable.  In particular, any two
 strong realisations of $A_*$ are canonically isomorphic.
\end{remark}

Our main aim is to prove that certain $R_*$-algebras are strongly
realisable, and to prove some more \emph{ad hoc} results for certain
algebras over $MU_*/2$.  

\begin{definition}
 A \emph{localised regular quotient (LRQ)} of $R_*$ is an algebra
 $A_*$ over $R_*$ that can be written in the form $A_*=(S^{-1}R_*)/I$,
 where $S$ is any set of (homogeneous) elements in $R_*$ and $I$ is an
 ideal which can be generated by a regular sequence.  We say that
 $A_*$ is a \emph{positive localised regular quotient (PLRQ)} if it
 can be written in the form $(S^{-1}R_*)/I$ as above, where $I$ can be
 generated by a regular sequence of elements of nonnegative degree.
\end{definition}

\begin{remark}
 If $A_*$ is an LRQ of $R_*$ and $B_*$ is an arbitrary $R_*$-algebra
 then $\CR_*(A_*,B_*)$ has at most one element.  Suppose that $A$ is a
 commutative ring object in $A\in\CR^{ec}$ with a given isomorphism
 $\pi_*(A)\simeq A_*$.  It follows that $A$ is a strong realisation of
 $A_*$ if and only if: whenever there is a map $A_*\xra{}\pi_*(B)$ of
 $R_*$-algebras, there is a unique map $A\xra{}B$ in $\CR^{ec}$.
\end{remark}

\begin{remark}\label{rem-localisation}
 Let $S$ be a set of homogeneous elements in $R_*$.  Using the results
 of~\cite[Section VIII.2]{ekmm:rma} one can construct a strictly
 commutative ring spectrum $S^{-1}R$ and a map $R\xra{}S^{-1}R$
 inducing an isomorphism $S^{-1}\pi_*(R)\xra{}\pi_*(S^{-1}R)$.
 Results of Wolbert show that $\DD_{S^{-1}R}$ is equivalent to the
 subcategory of $\DD_R$ consisting of objects $M$ such that each
 element of $S$ acts invertibly on $\pi_*(M)$.  Using this it is easy
 to check that any algebra over $S^{-1}R_*$ is strongly realisable
 over $R$ if and only if it is strongly realisable over $S^{-1}R$.
 For more discussion of this, see Section~\ref{sec-realise}.
\end{remark}

We start by stating a result for odd primes, which is relatively
easy. 
\begin{theorem}\label{thm-odd-general}
 If $A_*$ is an LRQ of $R_*$ and $2$ is a unit in $A_*$ then $A_*$ is
 strongly realisable.
\end{theorem}
This will be proved as Theorem~\ref{thm-odd-proof}.

Our main contribution is the extension to the case where $2$ is not
inverted.  Our results involve a certain ``commutativity obstruction''
$\cb(x)\in\pi_{2|x|+2}(R)/(2,x)$, which is defined in
Section~\ref{sec-prod-Rx}.  In Section~\ref{sec-pow-op}, we show that
when $d\geq 0$ this arises from a power operation
$\tP\:\pi_d(R)\xra{}\pi_{2d+2}(R)/2$.  This result was inspired by a
parallel result of Mironov in Baas-Sullivan theory~\cite{mi:mtc}.  The
restriction $d\geq 0$ is actually unneccessary but the argument for
the case $d<0$ is intricate and we have no applications so we have
omitted it.  In Section~\ref{sec-formal} we show how to compute this
power operation using formal group theory, at least in the case
$R=MU$.  The first steps in this direction were also taken by
Mironov~\cite{mi:mtc}, but our results are much more precise.

By Remark~\ref{rem-localisation} we also have a power operation
$\tP\:\pi_d(S^{-1}R)\xra{}\pi_{2d+2}(S^{-1}R)/2$.  This is in fact
determined algebraically by the power operation on $\pi_*R$, as we
will see in Section~\ref{sec-formal}.

Our result for the case where $A_*$ has no $2$-torsion is quite simple
and similar to the case where $2$ is inverted.
\begin{theorem}\label{thm-even-general}
 Let $A_*=(S^{-1}R_*)/I$ be a PLRQ of $R_*$ which has no $2$-torsion.
 Suppose also that $\tP(I)\subseteq S^{-1}R_*$ maps to $0$ in $A_*/2$.
 Then $A_*$ is strongly realisable.
\end{theorem}
This will be proved as Theorem~\ref{thm-even-proof}.

We next recall the definitions of some algebras over $MU_*$ which one
might hope to realise as spectra using the above results.  First, we
have the rings
\begin{align*}
 kU_*    & \eqdef \Zh[u] \qquad\qquad |u| = 2                      \\
 KU_*    & \eqdef \Zh[u^{\pm 1}]                                   \\
 H_*     & \eqdef \Zh_{\hphantom{p}}  \qquad\text{(in degree zero)}\\
 H\FF_*  & \eqdef \Fp  \qquad\text{(in degree zero)}.
\end{align*}
These are PLRQ's of $MU_*$ in well-known ways.  Next, we consider the
Brown-Peterson ring
\[ BP_* \eqdef \Zpl[v_k\st k>0] \qquad\qquad |v_k|=2(p^k-1). \]
We take $v_0=p$ as usual.  There is a unique $p$-typical formal group
law $F$ over this ring such that
\[ [p]_F(x) = \exp_F(px) +_F \sum^F_{k>0} v_k x^{p^k}. \]
(Thus, our $v_k$'s are Hazewinkel's generators rather than Araki's.)
We use this FGL to make $BP_*$ into an algebra over $MU_*$ in the
usual way.  We define
\begin{align*}
 P(n)_*  &\eqdef BP_*/(v_i\st i<n) = \Fp[v_j\st j\ge n]                \\
 B(n)_*  &\eqdef v_n^{-1}BP_*/(v_i\st i<n) = v_n^{-1}\Fp[v_j\st j\ge n]\\
 k(n)_*  &\eqdef BP_*/(v_i\st i\neq n) = \Fp[v_n]                      \\
 K(n)_*  &\eqdef v_n^{-1}BP_*/(v_i\st i\neq n) = \Fp[v_n^{\pm 1}]      \\
 \BP{n}_*&\eqdef BP_*/(v_i\st i>n) = \Zpl[v_1,\ldots,v_n]              \\
 E(n)_*  &\eqdef v_n^{-1}BP_*/(v_i\st i>n) =
              \Zpl[v_1,\ldots,v_{n-1},v_n^{\pm 1}]
\end{align*}
These are all PLRQ's of $BP_*$, and it is not hard to check that
$BP_*$ is a PLRQ of $MU_{(p)*}$, and thus that all the above rings are
PLRQ's of $MU_{(p)*}$.  

We also let $w_k\in\pi_{2(p^k-1)}MU$ denote the bordism class of a
smooth hypersurface $W_{p^k}$ of degree $p$ in $\cp^{p^k}$.  It is
well-known that $I_n=(w_i\st i<n)$ is the smallest ideal modulo which
the universal formal group law over $MU_*$ has height $n$, and that
the image of $I_n$ in $BP_*$ is the ideal $(v_i\st i<n)$.  In fact, we
have
\[ \sum_{m>0} [W_m] x^m dx = [p]_F(x) d\log_F(x) =
   [p]_F(x) \sum_{m\geq 0} [\cp^m] x^m dx.
\] 
Moreover, the sequence of $w_i$'s is regular, so that
$MU_*/I_n$ is a PLRQ of $MU_*$.

One can also define PLRQ's of $MU[\sixth]_*$ giving rise to various
versions of elliptic homology, but we refrain from giving details
here.  If we do not invert $6$ then the relevant rings seem not to be
LRQ's of $MU_*$.  If we take $R=MU^\wedge_p$ then we can make
$\Zh_p[v_n]$ into an LRQ of $R_*$ in such a way that the resulting
formal group law is of the (non-$p$-typical) type considered by Lubin
and Tate in algebraic number theory.  We can also take $R=L_{K(n)}MU$
and consider $\widehat{E(n)}_*$ as an LRQ of $R_*$ via the Ando
orientation~\cite{an:ifg} rather than the more usual $p$-typical one.
We leave the details of these applications to the reader.

The following proposition is immediate from
Theorem~\ref{thm-odd-general}.  
\begin{proposition}
 If $p>2$ and $R=MU$ or $R=MU_{(p)}$ then $kU_{(p)*}$, $KU_{(p)*}$,
 $H_{(p)*}$, $H\FF_{p*}$, $BP_*$, $P(n)_*$, $B(n)_*$, $k(n)_*$,
 $\BP{n}_*$, $E(n)_*$ and $MU_*/I_n$ are all strongly realisable.
\end{proposition}

After doing some computations with the power operation $\tP$, we will
also prove the following.
\begin{proposition}\label{prop-MU-omni}
 If $R=MU$ then $kU_*$, $KU_*$, $H_*$ and $H\FF_*$ are strongly
 realisable.  If $R=MU_{(2)}$ then $kU_{(2)*}$, $KU_{(2)*}$,
 $H_{(2)*}$ and $BP_*$ are strongly realisable.
\end{proposition}

The situation is less satisfactory for the rings $\BP{n}_*$ and
$E(n)_*$ at $p=2$.  For $n>1$, they cannot be realised as the homotopy
rings of commutative ring objects in $\DD$.  However, if we kill off a
slightly different sequence of elements instead of the sequence
$(v_{n+1},v_{n+2},\ldots)$, we get a quotient ring that is
realisable.  The resulting spectrum serves as a good substitute for
$\BP{n}$ in almost all arguments.  
\begin{proposition}\label{prop-BPn}
 If $R=MU_{(2)}$ and $n>0$, there is a quotient ring $\BP{n}'_*$ of
 $BP_*$ such that
 \begin{enumerate}
 \item The evident map
  \[ \Ztl[v_1,\ldots,v_n] \mra BP_* \era \BP{n}'_* \]
  is an isomorphism.
 \item $\BP{n}'_*$ is strongly realisable.
 \item We have $\BP{n}'_*/I_n=k(n)_*=BP_*/(v_i\st i\neq n)$ as
  $MU_*$-algebras.  
 \end{enumerate}
 Moreover, the ring $E(n)'_*=v_n^{-1}\BP{n}'_*$ is also strongly
 realisable.  If $n=1$ then we can take $\BP{1}'_*=\BP{1}_*$.
\end{proposition}
This is proved in Section~\ref{sec-mu}.

The situation for $MU_*/2$ and algebras over it is also more
complicated than for odd primes.
\begin{definition}
 Throughout this paper, we write $\tau$ for the twist map
 $X\Smash X\xra{}X\Smash X$, for any object $X$ for which this makes
 sense.  We say that a ring map $f\:A\xra{}B$ in $\CR$ is
 \emph{central} if
 \[ \phi\circ\tau\circ(f\Smash 1) = \phi\circ(f\Smash 1)\:
    A\Smash B\xra{} B,
 \]
 where $\phi\:B\Smash B\xra{}B$ is the product.  We say that $B$ is a
 \emph{central $A$-algebra} if there is a given central map
 $A\xra{}B$.
\end{definition}

\begin{theorem}\label{thm-MUIn}
 When $R=MU_{(2)}$, there is a ring $MU/I_n\in\CR$ with
 $\pi_*(MU/I_n)=MU_*/I_n$, and derivations
 $Q_i\:MU/I_n\xra{}\Sg^{2^{i+1}-1}MU/I_n$ for $0\leq i<n$.  If $\phi$
 is the product on $MU/I_n$ we have
 \[ \phi\circ\tau-\phi=w_n \phi\circ (Q_{n-1}\Smash Q_{n-1}). \]
\end{theorem}
This is proved in Section~\ref{sec-mu}.  There are actually many
non-isomorphic rings with these properties.  We will outline an
argument that specifies one of them unambiguously.

We get a sharper statement for algebras over $P(n)_*$.
\begin{theorem}\label{thm-Pn}
 When $R=MU_{(2)}$, there is a central $BP$-algebra
 $P(n)=BP\Smash MU/I_n\in\CR$ and an isomorphism $\pi_*P(n)=P(n)_*$.
 This has derivations  $Q_i\:P(n)\xra{}\Sg^{2^{i+1}-1}P(n)$ for
 $0\leq i<n$.  If $\phi$ is the product on $P(n)$ we have
 \[ \phi\circ\tau-\phi=v_n \phi\circ(Q_{n-1}\Smash Q_{n-1}). \]
 If $B$ is another central $BP$-algebra such that 
 \[ \pi_kB = \begin{cases}
      \{0,1\}     & \text{ if } k=0                     \\
      0           & \text{ if } 0<k<|v_n|               \\
      \{0,v_n\}   & \text{ if } k=|v_n|
    \end{cases}
 \]
 then either there is a unique map $P(n)\xra{}B$ of $BP$-algebras, or
 there is a unique map $P(n)\xra{}B^{\text{op}}$.  Analogous
 statements hold for $B(n)$, $k(n)$ and $K(n)$ with $BP$ replaced by
 $v_n^{-1}BP$, $\BP{n}'$ and $E(n)'$ respectively.
\end{theorem}
This is also proved in Section~\ref{sec-mu}.  Related results were
announced by W\"{u}rgler in~\cite{wu:crs}, but there appear to be some
problems with the line of argument used there.  A correct proof on
similar lines has recently been given by Nassau~\cite{na:spp,na:nks}.

\section{Products on $R/x$}\label{sec-prod-Rx}

Suppose that $x\in R_d$ is not a zero-divisor (so $d$ is even).  We
then have a cofibre sequence in the triangulated category $\DD$:
\[ \Sg^d R \xra{x} R \xra{\rho} R/x \xra{\bt} \Sg^{d+1}R. \]
Because $x$ is not a zero divisor, we have $\pi_*(R/x)=R_*/x$.  In
particular, $\pi_{d+1}R/x=0$ (because $d+1$ is odd), and thus
$\rho^*\:[R/x,R/x]\simeq[R,R/x]$.  It follows that $R/x$ is unique up
to unique isomorphism as an object under $R$.  

We next set up a theory of products on objects of the form $R/x$.
Apart from the fact that all such products are associative, our
results are at most minor sharpenings of the those 
in~\cite[Chapter V]{ekmm:rma}.

Observe that $\Rxx$ is a cell $R$-module with one $0$-cell, two
$(d+1)$-cells and one $(2d+2)$-cell.  We say that a map
$\phi\:\Rxx\xra{}R/x$ is a {\em product} if it agrees with $\rho$ on
the bottom cell, in other words
$\phi\circ(\rho\Smash\rho)=\rho\:R\xra{}R/x$.

The main result is as follows.
\begin{proposition}\label{prop-Rx}
 \begin{enumerate}
 \item All products are associative, and have $\rho$ as a two-sided
  unit.  
 \item The set of products on $R/x$ has a free transitive action of
  the group $R_{2d+2}/x$ (in particular, it is nonempty).
 \item There is a naturally defined element
  $\cb(x)\in\pi_{2d+2}(R)/(2,x)$ such that $R/x$ admits a commutative
  product if and only if $\cb(x)=0$.
 \item If so, the set of commutative products has a free transitive
  action of $\ann(2,R_{2d+2}/x)=\{y\in R_{2d+2}/x\st 2y=0\}$.
 \item If $d\geq 0$ there is a power operation
  $\tP\:R_d\xra{}R_{2d+2}/2$ such that $\cb(x)=\tP(x)\pmod{2,x}$ for
  all $x$.
 \end{enumerate}
\end{proposition}
\begin{proof}
 Part~(1) is proved as Lemma~\ref{lem-unital} and
 Proposition~\ref{prop-ass}.  In part~(2), the fact that products
 exist is~\cite[Theorem V.2.6]{ekmm:rma}; we also give a proof in
 Corollary~\ref{cor-products-exist}, which is slightly closer in
 spirit with our other proofs.  Parts~(3) and~(4) form
 Corollary~\ref{cor-comm}.  Part~(5) is explained in more detail and
 proved in Section~\ref{sec-pow-op}.
\end{proof}
From now on we will generally state our results in terms of $\tP(x)$
instead of $\cb(x)$, as that is the form in which the results are
actually applied.

\begin{lemma}\label{lem-x-zero}
 The map $x\:\Sg^dR/x\xra{}R/x$ is zero.
\end{lemma}
\begin{proof}
 Using the cofibration 
 \[ \Sg^dR\xra{\rho}\Sg^dR/x\xra{\bt}\Sg^{d+1}R \]
 and the fact that $\pi_{d+1}(R/x)=\pi_{d+1}(R)/x=0$, we find that
 $\rho^*\:[R/x,R/x]_d\xra{}[R,R/x]_d=\pi_d(R/x)$ is injective.  It is
 clear that $x$ gives zero on the right hand side, so it is zero on
 the left hand side as claimed.
\end{proof}
\begin{corollary}\label{cor-products-exist}
 There exist products on $R/x$.
\end{corollary}
\begin{proof}
 There is a cofibration $\Sg^dR/x\xra{x}R/x\xra{1\Smash\rho}\Rxx$.
 The lemma tells us that the first map is zero, so $1\Smash\rho$ is a
 split monomorphism, and any splitting is clearly a product.
\end{proof}

\begin{lemma}\label{lem-unital}
 If $\phi\:\Rxx\xra{}R/x$ is a product then $\rho$ is a two-sided unit
 for $\phi$, in the sense that
 \[ \phi\circ(\rho\Smash 1) = \phi\circ(1\Smash \rho) = 1
        \: R/x \xra{} R/x.
 \]
\end{lemma}
\begin{proof}
 By hypothesis, $\phi\circ(\rho\Smash 1)\:R/x\xra{}R/x$ is the
 identity on the bottom cell of $R/x$.  We observed earlier that
 $[R/x,R/x]\simeq[R,R/x]$, and it follows that
 $\phi\circ(\rho\Smash 1)=1$.  Similarly $\phi\circ(1\Smash \rho)=1$.
\end{proof}
\begin{remark}
 EKMM study products for which $\rho$ is a one-sided unit, and our
 definition of products is \emph{a priori} even weaker.  It follows
 from the lemma that EKMM's products are the same as ours and have
 $\rho$ as a two-sided unit.
\end{remark}

\begin{lemma}\label{lem-A-split}
 Let $A\in\DD$ be such that $x\:\Sg^dA\xra{}A$ is zero.  Then the
 diagram
 \[ R/x\Wedge R/x \xra{(\rho\Smash 1,1\Smash\rho)} \Rxx 
    \xra{\btt} \Sg^{2d+2}R
 \]
 induces a left-exact sequence
 \[ [\Sg^{2d+2}R,A] \mra [\Rxx,A] \xra{} [R/x\Wedge R/x,A]. \]
 Similarly, the diagram 
 \[ \Rxx\Wedge\Rxx\Wedge\Rxx 
    \xra{(\rho \Smash 1 \Smash 1,
          1 \Smash \rho \Smash 1,
          1 \Smash 1 \Smash \rho)}
    \Rxxx \xra{\bt^{(3)}} \Sg^{3d+3}R
 \]
 gives a left-exact sequence
 \[ [\Sg^{3d+3}R,A] \mra [\Rxxx,A] \xra{}
       [\Rxx\Wedge\Rxx\Wedge\Rxx,A].
 \]
\end{lemma}
\begin{proof}
 Consider the following diagram:
 \[ \diagram
   & R/x \dto^{\rho\Smash 1} \\
   R/x \dto_\bt \rto^{1\Smash\rho} &
   \Rxx \dto^{\bt\Smash 1} \drto^{\btt} \\
   \Sg^{d+1}R \rto_\rho &
   \Sg^{d+1}R/x \rto_\bt &
   \Sg^{2d+2} R
 \enddiagram \]
 We now apply the functor $[-,A]$ and make repeated use of the
 cofibration 
 \[ \Sg^d R \xra{x} R \xra{\rho} R/x \xra{\bt} \Sg^{d+1} R. \]
 The conclusion is that all maps involving $\bt$ become monomorphisms,
 all maps involving $\rho$ become epimorphisms, and the bottom row and
 the middle column become short exact.  The first claim follows by
 diagram chasing.  For the second claim, consider the diagram
 \[ \diagram
   & \Rxx\Wedge\Rxx \dto^{(\rho\Smash 1\Smash 1,1\Smash\rho\Smash 1)}\\
   \Rxx \dto_{\btt} \rto^{1\Smash 1\Smash\rho} &
   \Rxxx \dto^{\btt\Smash 1} \drto^{\btt\Smash\bt} \\
   \Sg^{2d+2}R \rto_\rho &
   \Sg^{2d+2}R/x \rto_\bt &
   \Sg^{3d+3} R
 \enddiagram \]
 We apply the same logic as before, using the first claim (with $A$
 replaced by $F(R/x,A)$) to see that the middle column becomes left
 exact. 
\end{proof}

We next determine how many different products there are on $R/x$.
\begin{lemma}\label{lem-uni-obs}
 If $\phi$ is a product on $R/x$ and
 $u\in\pi_{2d+2}(R)/x=[\Sg^{2d+2}R,R/x]$ then
 $\phi'=\phi+u\circ(\btt)$ is another product.  Moreover, this
 construction gives a free transitive action of $\pi_{2d+2}(R)/x$ on
 the set of all products.
\end{lemma}
\begin{proof}
 Let $P$ be the set of products.  As $(\btt)\circ(\rho\Smash\rho)=0$,
 it is clear that the above construction gives an action of
 $\pi_{2d+2}(R)/x$ on $P$.  Now suppose that $\phi,\phi'\in P$.  We
 need to show that there is a unique $u\:\Sg^{2d+2}R\xra{}R/x$ such
 that $\phi'=\phi+u\circ(\btt)$.  Using the unital properties of
 $\phi$ and $\phi'$ given by Lemma~\ref{lem-unital}, we see that
 \[ (\phi'-\phi)\circ(\rho\Smash 1)=(\phi'-\phi)\circ(1\Smash\rho)=0.
 \]
 Because of Lemma~\ref{lem-x-zero}, we can apply
 Lemma~\ref{lem-A-split} to see that $\phi'-\phi=u\circ(\btt)$ for a
 unique element $u$, as claimed.
\end{proof}

\begin{proposition}\label{prop-ass}
 Any  product on $R/x$ is associative.
\end{proposition}
\begin{proof}
 Let $\phi$ be a product, and write 
 \[ \dl\eqdef\phi\circ(\phi\Smash 1-1\Smash\phi)\:\Rxxx \xra{} R/x,
 \]
 so the claim is that $\dl$ is nullhomotopic.  Using the unital
 properties of $\phi$ we see that
 \[ \dl\circ(\rho\Smash 1\Smash 1) = 
    \dl\circ(1 \Smash\rho\Smash 1) = 
    \dl\circ(1 \Smash 1\Smash\rho) = 0.
 \]
 Using Lemma~\ref{lem-A-split}, we conclude that
 $\dl=u\circ(\bt\Smash\btt)$ for a unique element
 $u\in[\Sg^{3d+3}R,R/x]=\pi_{3d+3}(R)/x=0$ (because $3d+3$ is odd).
 Thus $\dl=0$ as claimed.
\end{proof}
\begin{remark}
 The corresponding result in Baas-Sullivan theory was already known
 (this is proved in~\cite{mi:ems} in a form which is valid when $R_*$
 need not be concentrated in even degrees, for example for $R=MSp$).
\end{remark}
\begin{remark}\label{rem-ass}
 The EKMM approach to associativity is essentially as follows.  They
 note that $R/x$ has cells of dimension $0$ and $d+1$, so
 $(R/x)^{(3)}$ has cells in dimensions $0$, $d+1$, $2d+2$ and $3d+3$.
 The map $\dl$ vanishes on the zero-cell and
 $\pi_{d+1}(R/x)=\pi_{3d+3}(R/x)=0$ so the only obstruction to
 concluding that $\dl=0$ lies in $\pi_{2d+2}(R/x)$.  EKMM work only
 with LRQ's that are concentrated in degrees divisible by $4$, so the
 obstruction goes away.  We instead use Lemma~\ref{lem-A-split} to
 analyse the attaching maps in $(R/x)^{(3)}$; implicitly, we show that
 the obstruction is divisible by $x$ and thus is zero.
\end{remark}

We now discuss commutativity.
\begin{lemma}\label{lem-comm}
 There is a natural map $c$ from the set of products to
 $\pi_{2d+2}R/x$ such that $c(\phi)=0$ if and only if $\phi$ is
 commutative.  Moreover,
 \[ c(\phi + u\circ(\btt)) = c(\phi) - 2u.  \]
\end{lemma}
\begin{proof}
 Let $\tau\:\Rxx\xra{}\Rxx$ be the twist map.  Clearly, if $\phi$ is a
 product then so is $\phi\circ\tau$.  Thus, there is a unique element
 $v\in\pi_{2d+2}R/x$ such that
 \[ \phi\circ\tau = \phi + v\circ(\btt). \]
 We define $c(\phi)\eqdef v$.  Next, recall that the twist map on
 $\Sg^{2d+2}R=\Sg^{d+1}R\Smash\Sg^{d+1}R$ is homotopic to $(-1)$,
 because $d+1$ is odd.  It follows by naturality that
 $(\btt)\circ\tau=\tau\circ(\btt)=-\btt$.  Consider a second product
 $\phi'=\phi+u\circ(\btt)$.  We now see that
 \[ \phi'\circ\tau = \phi + v\circ(\btt) - u\circ(\btt) =
        \phi' + (v-2u)\circ(\btt).
 \]
 Thus $c(\phi')=c(\phi)-2u$ as claimed.  
\end{proof}

\begin{corollary}\label{cor-comm}
 There is a naturally defined element $\cb(x)\in\pi_{2d+2}(R)/(2,x)$
 such that $R/x$ admits a commutative product if and only if
 $\cb(x)=0$.  If so, the set of commutative products has a free
 transitive action of the group
 $\ann(2,\pi_{2d+2}(R)/x)\eqdef\{y\in\pi_{2d+2}(R)/x\st 2y=0\}$.  In
 particular, if $\pi_*(R)/x$ has no $2$-torsion then there is a unique
 commutative product.
\end{corollary}
\begin{proof}
 We choose a product $\phi$ on $R/x$ and define
 $\cb(x)\eqdef c(\phi)\pmod{2}$.  This is well-defined, by the lemma.
 If $\cb(x)\neq 0$ then $c(\phi')\neq 0$ for all $\phi'$, so there is
 no commutative product.  If $\cb(x)=0$ then $c(\phi)=2w$, say, so
 that $\phi'=\phi+w\circ(\btt)$ is a commutative product.  In this
 case, the commutative products are precisely the products of the form
 $\phi'+z\circ(\btt)$ where $2z=0$, so they have a free transitive
 action of $\ann(2,\pi_{2d+2}(R)/x)$.
\end{proof}

Next, we consider the Bockstein operation:
\[ \btb \eqdef \rho\bt\: R/x \xra{} \Sg^{d+1} R/x. \]

\begin{definition}
 Let $A\in\CR$ be a ring, with product $\phi\:A\Smash A\xra{}A$.  We
 say that a map $Q\:A\xra{}\Sg^kA$ is a \emph{derivation} if we have
 \[ Q\circ\phi = \phi\circ(Q\Smash 1 + 1\Smash Q)\: A^{(2)}\xra{}A.
 \]
\end{definition}

\begin{proposition}\label{prop-der}
 The map $\btb$ is a derivation with respect to any product $\phi$ on
 $R/x$.
\end{proposition}
\begin{proof}
 Write $\dl\eqdef\btb\circ\phi-\phi\circ(\btb\Smash 1+1\Smash\btb)$,
 so the claim is that $\dl=0$.  It is easy to see that
 $\dl\circ(\rho\Smash 1)=\dl\circ(1\Smash\rho)=0$, so by
 Lemma~\ref{lem-A-split} we see that $\dl$ factors through a unique
 map $\Sg^{2d+2}R\xra{}\Sg^{d+1}R/x$.  This is an element of
 $\pi_{d+1}(R)/x$, which is zero because $d+1$ is odd.
\end{proof}

We end this section by analysing maps out of the rings $R/x$.
\begin{proposition}\label{prop-maps-Rx}
 Let $A\in\CR^e$ be an even ring.  If $x$ maps to zero in $\pi_*A$
 then there is precisely one unital map $f\:R/x\xra{}A$, and otherwise
 there are no such maps.  If $f$ exists and $\phi$ is a product on
 $R/x$, then there is a naturally defined element
 $d_A(\phi)\in\pi_{2d+2}(A)$ such that
 \begin{itemize}
 \item[(a)] $d_A(\phi)=0$ if and only if $f$ is a ring map with respect to
  $\phi$.  
 \item[(b)] $d_A(\phi+u\circ(\btt))=d_A(\phi)+u$.
 \item[(c)] If $A$ is commutative then
  $2d_A(\phi)=c(\phi)\in\pi_{2d+2}A$.
 \end{itemize}
\end{proposition}
\begin{proof}
 The statement about the existence and uniqueness of $f$ follows
 immediately from the cofibration
 $\Sg^dR\xra{x}R\xra{\rho}R/x\xra{\bt}\Sg^{d+1}R$, and the fact that
 $\pi_{d+1}A=0$.  Suppose that $f$ exists; it follows easily using the
 product structure on $A$ that $x\:\Sg^dA\xra{}A$ is zero.  Now let
 $\psi$ be the given product on $A$, and let $\phi$ be a product on
 $R/x$.  Consider the map
 \[ \dl\eqdef\psi\circ(f\Smash f)-f\circ\phi\:\Rxx\xra{}A. \]
 By the usual argument, we have $\dl=v\circ(\btt)$ for a unique map
 $v\:\Sg^{2d+2}R\xra{}A$.  We define
 $d_A(\phi)\eqdef v\in\pi_{2d+2}A$.  It is obvious that this vanishes
 if and only if $f$ is a ring map, and that
 $d_A(\phi+u\circ(\btt))=d_A(\phi)+u$.

 Now suppose that $A$ is commutative, so $\psi=\psi\circ\tau$.  On the
 one hand, using the fact that $(\btt)\circ\tau=-\btt$ we see that
 $\dl-\dl\circ\tau=2d_A(\phi)\circ(\btt)$.  On the other hand, from
 the definition of $\dl$ and the fact that $\psi\circ\tau=\psi$, we
 see that
 \[ \dl-\dl\circ\tau=f\circ(\phi-\phi\circ\tau)=c(\phi)\circ(\btt).
 \]
 Because $(\btt)^*\:\pi_{2d+2}A\xra{}[\Rxx,A]$ is a split
 monomorphism, we conclude that $2d_A(\phi)=c(\phi)$ in
 $\pi_{2d+2}A$. 
\end{proof}

\section{Strong realisations}\label{sec-realise}

In this section we assemble the products which we have constructed on
the $R$-modules $R/x$ to get products on more general $R_*$-algebras.
We will work entirely in the derived category $\DD$, rather than the
underlying geometric category.  All the main ideas in this section
come from~\cite[Chapter V]{ekmm:rma}.

We start with some generally nonsensical preliminaries.  
\begin{definition}
 Given a diagram $A\xra{f}C\xla{g}B$ in $\CR$, we say that $f$
 commutes with $g$ if and only if we have 
 \[ \phi_C\circ(f\Smash g) = \phi_C\circ\tau\circ(f\Smash g) 
    \: A\Smash B \xra{} C . 
 \]
 Note that this can be false when $f=g$; in particular $A$ is
 commutative if and only if $1_A$ commutes with itself.
\end{definition}

The next three lemmas become trivial if we replace $\DD$ by the
category of modules over a commutative ring, and the smash product by
the tensor product.  The proofs in that context can easily be made
diagrammatic and thus carried over to $\DD$.
\begin{lemma}\label{lem-ring-smash}
 If $A$ and $B$ are rings in $\CR$, then there is a unique ring
 structure on $A\Smash B$ such that the evident maps
 $A\xra{i}A\Smash B\xla{j}B$ are commuting ring maps.  Moreover, with
 this product, $(i,j)$ is the universal example of a commuting pair of
 maps out of $A$ and $B$.  \qed
\end{lemma}

\begin{lemma}\label{lem-comm-test}
 A map $f\:A\Smash B\xra{}C$ commutes with itself if and only if
 $f\circ i$ commutes with itself and $f\circ j$ commutes with itself.
 In particular, $A\Smash B$ is commutative if and only if $i$ and $j$
 commute with themselves.  \qed
\end{lemma}

\begin{lemma}\label{lem-comm-smash}
 If $A$ and $B$ are commutative, then so is $A\Smash B$, and it is the
 coproduct of $A$ and $B$ in $\CR^c$.  \qed
\end{lemma}

\begin{corollary}\label{cor-tensor}
 Suppose that
 \begin{itemize}
  \item $A$ and $B$ are strong realisations of $A_*$ and $B_*$.
  \item The ring $A_*\ot_{R_*}B_*$ has no $2$-torsion.
  \item The natural map $A_*\ot_{R_*}B_*\xra{}\pi_*(A\Smash B)$ is an
   isomorphism. 
 \end{itemize}
 Then $A\Smash B$ is a strong realisation of $A_*\ot_{R_*}B_*$. \qed
\end{corollary}

We next consider the problem of realising $S^{-1}R_*$, where $S$ is a
set  of homogeneous elements of $R_*$.  If $S$ is countable then we
can construct an object $S^{-1}R\in\DD$ by the method
of~\cite[Section V.2]{ekmm:rma}; this has
$\pi_*(S^{-1}R)=S^{-1}\pi_*(R)$.  If we want to allow $S$ to be
uncountable then it seems easiest to construct $S^{-1}R$ as the finite
localisation of $R$ away from the $R$-modules $\{R/x\st s\in S\}$;
see~\cite{mi:fl} or~\cite[Theorem 3.3.7]{hopast:ash}.  In either case,
we note that $S^{-1}R$ is the Bousfield localisation of $R$ in $\DD$
with respect to $S^{-1}R$.  We may thus 
use~\cite[Section VIII.2]{ekmm:rma} to construct a model of $S^{-1}R$
which is a strictly commutative algebra over $R$ in the underlying
topological category of spectra.  The localisation functor involved
here is smashing, so results of Wolbert~\cite{wo:cmk} 
\cite[Section VIII.3]{ekmm:rma} imply that $\DD_{S^{-1}R}$ is
equivalent to the full subcategory of $\DD_R$ consisting of
$R$-modules $M$ for which $\pi_*(M)$ is a module over $S^{-1}R_*$.  
This makes the following result immediate.

\begin{proposition}\label{prop-S-inv}
 Let $S$ be a set of homogeneous elements of $R_*$, and let $A_*$ be
 an algebra over $S^{-1}R_*$.  Then $A_*$ is strongly realisable over
 $R$ if and only if it is strongly realisable over $S^{-1}R$. \qed
\end{proposition}
This allows us to reduce everything to the case $S=\emptyset$.

Now consider a sequence $(x_i)$ in $R_*$, with products $\phi_i$ on
$R/x_i$.  Write $A_i=R/x_1\Smash\ldots\Smash R/x_i$, and make this
into a ring as in Lemma~\ref{lem-ring-smash}.  There are evident maps
$A_i\xra{}A_{i+1}$, so we can form the telescope $A=\hcolim_iA_i$.  

\begin{lemma}\label{lem-lim-one}
 If $M\in\DD$ and $I=(x_1,x_2,\ldots)\leq R_*$ acts trivially on
 $M$ and $r\geq 0$ then $[A^{(r)},M]=\invlim_i[A^{(r)}_i,M]$.
\end{lemma}
\begin{proof}
 This will follow immediately from the Milnor sequence if we can show
 that $\invlim_i^1[A^{(r)}_i,M]_*=0$.  For this, it suffices to show
 that the map $\rho^*\:[B\Smash R/x_i,M]\xra{}[B,M]$ is surjective for
 all $B$.  This follows from the cofibration
 $\Sg^{|x_i|}B\xra{x_i}B\xra{}B\Smash R/x_i$ and the fact that $x_i$
 acts trivially on $M$.
\end{proof}

\begin{proposition}\label{prop-coprod}
 Let $(x_i)$ be a sequence in $R_*$, and $\phi_i$ a product on $R/x_i$
 for each $i$.  Let $A$ be the homotopy colimit of the rings
 $A_i=R/x_1\Smash\ldots\Smash R/x_i$, and let $f_i\:R/x_i\xra{}A$ be
 the evident map.  Then there is a unique associative and unital
 product on $A$ such that maps $f_i$ are ring maps, and $f_i$ commutes
 with $f_j$ when $i\neq j$.  This product is commutative if and only
 if each $f_i$ commutes with itself.  Ring maps from $A$ to any ring
 $B$ biject with systems of ring maps $g_i\:R/x_i\xra{}B$ such that
 $g_i$ commutes with $g_j$ for all $i\neq j$.
\end{proposition}
\begin{proof}
 Because $R/x_i$ admits a product, we know that
 $x_i$ acts trivially on $R/x_i$.  Because $A$ has the form
 $R/x_i\Smash B$, we see that $x_i$ acts trivially on $A$.  Thus $I$
 acts trivially on $A$, and Lemma~\ref{lem-lim-one} assures us that
 $[A^{(r)},A]=\invlim_i[A_i^{(r)},A]$.  

 Let $\psi_i$ be the product on $A_i$.  By the above, there is a
 unique map $\psi\:A\Smash A\xra{}A$ which is compatible with the maps
 $\psi_i$.  It is easy to check that this is an associative and unital
 product, and that it is the only one for which the $f_i$ are
 commuting ring maps.  It is also easy to check that $\psi$ is
 commutative if and only if each of the maps $A_i\xra{}A$ commutes
 with itself, if and only if each $f_i$ commutes with itself.

 Now let $B$ be any ring in $\CR$.  We may assume that each $x_i$ maps
 to zero in $\pi_*(B)$, for otherwise the claimed bijection is between
 empty sets.  As $B$ is a ring, this means that each $x_i$ acts
 trivially on $B$, so that $[A^{(r)},B]=\invlim_i[A_i^{(r)},B]$.  We
 see from Lemma~\ref{lem-ring-smash} that ring maps from $A_i$ to $B$
 biject with systems of ring maps $g_j\:R/x_j\xra{}B$ for $j<i$ such
 that $g_j$ commutes with $g_k$ for $j\neq k$.  The claimed
 description of ring maps $A\xra{}B$ follows easily.
\end{proof}

\begin{corollary}\label{cor-coprod}
 If each $R/x_i$ is commutative, then $A$ is the coproduct of the
 $R/x_i$ in $\CR^c$. \qed
\end{corollary}

\begin{remark}
 If the sequence $(x_i)$ is regular, then it is easy to see that
 $\pi_*(A)=R_*/(x_1,x_2,\ldots)$.  Note also that ring maps out of
 $R/x$ were analysed in Proposition~\ref{prop-maps-Rx}.
\end{remark}

We now  restate and prove Theorems~\ref{thm-odd-general}
and~\ref{thm-even-general}.  Of course, the former is a special case
of the latter, but it seems clearest to prove
Theorem~\ref{thm-odd-general} first and then explain the improvements
necessary for Theorem~\ref{thm-even-general}.

\begin{theorem}\label{thm-odd-proof}
 If $A_*$ is an LRQ of $R_*$ and $2$ is a unit in $A_*$ then $A_*$ is
 strongly realisable.
\end{theorem}

\begin{proof}
 We can use Proposition~\ref{prop-S-inv} to reduce to the case where
 $A_*=R_*/I$ where $2$ is invertible in $R_*$ and $I$ is generated by
 a regular sequence $(x_1,x_2,\ldots)$.  We know from
 Proposition~\ref{prop-Rx} that there is a unique commutative product
 $\phi_i$ on $R/x_i$.  If $C\in\CR^{ec}$ and $x_i=0$ in $\pi_*(C)$
 then in the notation of Proposition~\ref{prop-maps-Rx} we have
 $2d_C(\phi_i)=0$ and thus $d_C(\phi_i)=0$, so the unique unital map
 $R/x_i\xra{}C$ is a ring map.  It follows that $R/x_i$ is a strong
 realisation of $R_*/x_i$, and thus that
 $A_i=R/x_1\Smash\ldots\Smash R/x_i$ is a strong realisation of
 $R_*/(x_1,\ldots,x_i)$.  Using Proposition~\ref{prop-coprod}, we get
 a ring $A$ which is a strong realisation of $R_*/I$.
\end{proof}

We next address the case where $2$ is not a zero-divisor, but is not
invertible either.
\begin{theorem}\label{thm-even-proof}
 Let $A_*=(S^{-1}R_*)/I$ be a PLRQ of $R_*$ which has no $2$-torsion.
 Suppose also that $\tP(I)\subseteq S^{-1}R_*$ maps to $0$ in $A_*/2$,
 where $\tP\:R_d\xra{}R_{2d+2}/2$ is the power operation defined in
 Section~\ref{sec-pow-op}.  Then $A_*$ is strongly realisable.
\end{theorem}
\begin{proof}
 After using Proposition~\ref{prop-S-inv}, we may assume that
 $S=\emptyset$.  Choose a regular sequence $(x_i)$ generating $I$.  As
 $\cb(x_i)=\tP(x_i)\in I\pmod{2}$, we can choose a product $\phi_i$ on
 $R/x_i$ such that $c(\phi_i)\in I$.  We let $A$ be the ``infinite
 smash product'' of the $R/x_i$, as in Proposition~\ref{prop-coprod},
 so that $\pi_*(A)=A_*$.  Because $c(\phi_i)$ maps to zero in
 $\pi_*(A)$, we see easily that the map $R/x_i\xra{}A$ commutes with
 itself.  By Proposition~\ref{prop-coprod}, we conclude that $A$ is
 commutative.

 Let $B\in\CR^{ec}$ be an even commutative ring, and that $\pi_*(B)$ has no
 $2$-torsion.  The claim is that $\CR(A,B)=\CR_*(A_*,\pi_*(B))$.  The
 right hand side has at most one element, and if it is empty, then the
 left hand side is also.  Thus, we may assume that there is a map
 $A_*\xra{}\pi_*(B)$ of $R_*$-algebras, and we need to show that there
 is a unique ring map $A\xra{}B$.  

 By Proposition~\ref{prop-coprod}, we know that ring maps $A\xra{}B$
 biject with systems of ring maps $R/x_i\xra{}B$ (which automatically
 commute as $B$ is commutative).  There is a unique unital map
 $f\:R/x_i\xra{}B$, and Proposition~\ref{prop-maps-Rx} tells us that
 the obstruction to $f$ being a homomorphism satisfies
 $2d_B(\phi_i)=c(\phi_i)=0\in \pi_*(B)$.  Because $\pi_*(B)$ has no
 $2$-torsion, we have $d_B(\phi_i)=0$, so there is a unique ring map
 $R/x_i\xra{}B$, and thus a unique ring map $A\xra{}B$ as required.
\end{proof}

The following result is also useful.
\begin{proposition}\label{prop-free}
 Let $A_*$ be a strongly realisable $R_*$-algebra, and let
 $A_*\xra{}B_*$ be a map of $R_*$-algebras that makes $B_*$ into a
 free module over $A_*$.  Then $B_*$ is strongly realisable.
\end{proposition}
\begin{proof}
 First, observe that if $F$ and $M$ are $A$-modules, there is a
 natural map 
 \[ \Hom_A(F,M) \xra{} \Hom_{A_*}(F_*,M_*), \]
 which is an isomorphism if $F$ is a wedge of suspensions of $A$ (in
 other words, a free $A$-module).

 Choose a homogeneous basis $\{e_i\}$ for $B_*$ over $A_*$, where
 $e_i$ has degree $d_i$.  Define $B\eqdef\bigvee_i\Sg^{d_i}A$, so that
 $B$ is a free $A$-module with a given isomorphism $\pi_*B\simeq B_*$
 of $A_*$-modules.  Define $B_0\eqdef A$ and $B_1\eqdef B$ and
 \begin{align*}
   B_2 &\eqdef \bigvee_{i,j}\Sg^{d_i+d_j}A      \\
   B_3 &\eqdef \bigvee_{i,j,k}\Sg^{d_i+d_j+d_k}A.
 \end{align*}
 The product map $\mu\:A\Smash A\xra{}A$ gives rise to evident maps
 $\phi_k\:B^{(k)}\xra{}B_k$ which in turn give isomorphisms
 $B_*^{\ot_{A_*}k}=\pi_*B_k$ of $A_*$-modules.  The multiplication
 map $B_*\ot_{A_*}B_*\xra{}B_*$ corresponds under the isomorphism
 $\Hom_A(B_2,B)=\Hom_{A_*}(\pi_*B_2,B_*)$ to a map $B_2\xra{}B$.
 After composing this with $\phi_2$, we get a product map
 $\mu_B\:B\Smash B\xra{}B$.  A similar procedure gives a unit map
 $A\xra{}B$.

 We next prove that this product is associative.  Each of the two
 associated products $B^{(3)}\xra{}B$ factors as $\phi_3$ followed by
 a map $B_3\xra{}B$, corresponding to a map
 $B_*^{\ot_{A_*}3}\xra{}A_*$.  The two maps
 $B_*^{\ot_{A_*}3}\xra{}A_*$ in question are just the two possible
 associated products, which are the same because $B_*$ is associative.
 It follows that $B$ is associative.  Similar arguments show that $B$
 is commutative and unital.

 Now consider an object $C\in\CR$ equipped with a map $B_*\xra{}C_*$
 (and thus a map $A_*\xra{}C_*$).  As $A$ is a strong realisation of
 $A_*$, there is a unique map $A\xra{}C$ compatible with the map
 $A_*\xra{}C_*$.  This makes $C$ into an $A$-module, and thus gives an
 isomorphism $\Hom_A(B,C)=\Hom_{A_*}(B_*,C_*)$.  There is thus a
 unique $A$-module map $B\xra{}C$ inducing the given map
 $B_*\xra{}C_*$.  It follows easily that $B_*$ is a strong realisation
 of $B_*$.
\end{proof}

We will need to consider certain $R_*$-algebras that are not strongly
realisable.  The following result assures us that weaker kinds of
realisation are not completely uncontrolled.
\begin{proposition}\label{prop-weak-unique}
 Let $A_*$ be an LRQ of $R_*$, and let $B,C\in\CR^e$ be rings (not
 necessarily commutative) such that $\pi_*(B)=A_*=\pi_*(C)$.  Then
 there is an isomorphism $f\:B\xra{}C$ (not necessarily a ring
 map) that is compatible with the unit maps $B\xla{}R\xra{}C$.
\end{proposition}
\begin{proof}
 We may as usual assume that $S=\emptyset$, and write
 $I=(x_1,x_2,\ldots)$.  Let $A$ be the infinite smash product of the
 $R/x_i$'s, so that $\pi_*(A)=A_*$.  It will be enough to show that
 there is a unital isomorphism $A\xra{}B$.  Moreover, any unital map
 $A\xra{}B$ is automatically an isomorphism, just by looking at the
 homotopy groups.

 There is a unique unital map $f_i\:R/x_i\xra{}B$.  Write
 $A_i=R/x_1\Smash\ldots\Smash R/x_i$, and let $g_i$ be the map
 \[ A_i \xra{f_1\Smash\ldots\Smash f_i} B^{(i)} \xra{} B, \]
 where the second maps is the product.  Because $B$ is a ring and each
 $x_i$ goes to zero in $\pi_*(B)$, we can apply
 Lemma~\ref{lem-lim-one} to get a unital map $g\:A\xra{}B$ as
 required.
\end{proof}

We conclude this section by investigating $R$-module maps $A\xra{}A$
for various $R$-algebras $A\in\CR$.  
\begin{proposition}\label{prop-self-maps}
 Let $\{x_1,x_2,\ldots\}$ be a regular sequence in $R_*$, let $\phi_i$
 be a product on $R/x_i$, and let $A$ be the infinite smash product of
 the rings $R/x_i$.  Let $Q_i\:A\xra{}\Sg^{|x_i|+1}A$ be obtained by
 smashing the Bockstein map
 $\btb_{x_i}\:R/x_i\xra{}\Sg^{|x_i|+1}R/x_i$ with the identity map on
 all the other $R/x_j$'s.  Then $\DD(A,A)_*$ is isomorphic as an
 algebra over $A_*$ to the completed exterior algebra on the elements
 $Q_i$.
\end{proposition}
\begin{proof}
 It is not hard to see that $Q_iQ_j=-Q_jQ_i$, with a sign coming from
 an implicit permutation of suspension coordinates.  We also have
 $\btb_i^2=0$ and thus $Q_i^2=0$.  Given any finite subset
 $S=\{i_1<\ldots<i_n\}$ of the positive integers, we define 
 \[ Q_S\eqdef Q_{i_1} Q_{i_2} \ldots Q_{i_n} \: A \xra{}\Sg^{d_S} A, 
 \]
 where $d_S=\sum_j(|x_{i_j}|+1)$.  The claim is that one can make
 sense of homogeneous infinite sums of the form $\sum_S a_S Q_S$ with
 $a_S\in A_*$, and that any graded map $A\xra{}A$ of $R$-modules is
 uniquely of that form.

 Write $A_n=R/x_1\Smash\ldots\Smash R/x_n$, and let $i_n\:A_n\xra{}A$
 be the evident map.  It is easy to check that $Q_S\circ i_n=0$ if
 $\max(S)>n$, and a simple induction shows that $\DD(A_n,A)_*$ is a
 free module over $A_*$ generated by the maps $Q_S\circ i_n$ for which
 $\max(S)\leq n$.  Moreover, Lemma~\ref{lem-lim-one} implies that
 $\DD(A,A)_*=\invlim_n\DD(A_n,A)_*$.  The claim follows easily.  
\end{proof}

The above result relies more heavily than one would like on the choice
of a regular sequence generating the ideal $\ker(R_*\xra{}A_*)$.  We
will use the following construction to make things more canonical.
\begin{construction}\label{cons-dQ}
 Let $A\in\CR^e$ be an even ring, with unit $\eta\:R\xra{}A$, and let
 $I$ be the kernel of $\eta_*\:R_*\xra{}A_*$.  Given a derivation
 $Q\:A\xra{}\Sg^kA$, we define a function $d(Q)\:I\xra{}A_*$ as
 follows.  Given $x\in I$, we have a cofibration
 \[ \Sg^d R \xra{x} R \xra{\rho_x} R/x \xra{\bt_x} \Sg^{d+1}R \]
 as usual.  Here $x$ may be a zero-divisor in $R_*$, so we need not
 have $\pi_*(R/x)=\pi_*(R)/x$.  Nonetheless, we see easily that there
 is a unique map $f_x\:R/x\xra{}A$ such that $f_x\circ\rho_x=\eta$.
 As $Q$ is a derivation, one checks easily that $Q\circ\eta=0$, so
 $(Q\circ f_x)\circ\rho_x=0$, so $Q\circ f_x=y\circ\bt_x$ for some
 $y\:\Sg^{d+1}R\xra{}\Sg^kA$.  Because $x$ acts as zero on $A$, we see
 that $y$ is unique.  We can thus define
 $d(Q)(x)\eqdef y\in\pi_{d+1-k}A$.
\end{construction}

\begin{proposition}\label{prop-dQ}
 Let $A\in\CR^e$ be such that $\pi_*(A)=R_*/I$, where $I$ can be
 generated by a regular sequence.  Let $\Der(A)$ be the set of
 derivations $A\xra{}A$.  Then Construction~\ref{cons-dQ} gives rise
 to a natural monomorphism $d\:\Der(A)\xra{}\Hom_{R_*}(I/I^2,A_*)$
 (with degrees shifted by one).  
\end{proposition}
\begin{proof}
 Choose a regular sequence $\{x_1,x_2,\ldots\}$ generating $I$.
 Write $A_n=R/x_1\Smash\ldots\Smash R/x_n$, and let $j_n$ be the map 
 \[ R/x_1\Smash\ldots\Smash R/x_n \xra{f_{x_1}\Smash\ldots\Smash f_{x_n}}
    A^{(n)} \xra{\text{product}} A.
 \]
 It is easy to see that $A$ is the homotopy colimit of the objects
 $A_n$ (although there may not be a ring structure on $A_n$ for which
 $j_n$ is a homomorphism).  We also write $A_{n,i}$ for the smash
 product of the $R/x_j$ for which $j\leq n$ and $j\neq i$, and
 $j_{n,i}$ for the evident map $A_{n,i}\xra{}A_n\xra{j_n}A$. 

 Consider a derivation $Q\:A\xra{}\Sg^kA$, and write $b_i=d(Q)(x_i)$.
 Because $Q$ is a derivation, we see that $Q\circ j_n$ is a sum of $n$
 terms, of which the $i$'th is $b_i$ times the composite
 \[ A_n = A_{n,i}\Smash R/x_i \xra{1\Smash\bt_{x_i}} 
    \Sg^{|x_i|+1} A_{n,i} \xra{j_{n,i}} \Sg^{|x_i|+1}A.
 \]

 Now consider an element $x=\sum_{i=1}^n a_ix_i$ of $I$.  It is easy
 to see that there is a unique unital map $f'_x\:R/x\xra{}A_n$, and
 that $j_n\circ f'_x=f_x$.  Now consider the following diagram.
 \[ \diagram
  \Sg^d R \rto^x \dto_{a_i} &
  R \rto ^{\rho_x} \dto_1 &
  R/x \rto^{\bt_x} \dto^{f'_x} &
  \Sg^{d+1}R \dto^{a_i} \\
  \Sg^{|x_i|} A_{n,i} \rto_{x_i} &
  A_{n,i} \rto_{1\Smash\rho_{x_i}} &
  A_n \rto_{1\Smash\bt_{x_i}} &
  \Sg^{|x_i|+1}A_{n,i}
 \enddiagram \]
 The left hand square commutes because the terms $a_jx_j$ for
 $j\neq i$ become zero in $\pi_*(A_{n,i})$.  It follows that there
 exists a map $R/x\xra{}A_n$ making the whole diagram commute.
 However, $f'_x$ is the \emph{unique} map making the middle square
 commute, so the whole diagram commutes as drawn.  Thus
 $j_{n,i} \circ (\bt_{x_i}\Smash 1) \circ f'_x = a_i\circ\bt_x$
 (thinking of $a_i$ as an element of $\pi_*(A)$).  As
 $Q\circ j_n=\sum_i b_i.(j_{n,i}\circ(\bt_{x_i}\Smash 1))$, we
 conclude that
 $Q\circ f_x=Q\circ j_n\circ f'_x=(\sum_i a_ib_i)\circ\bt_x$.  Thus
 $d(Q)(x)=\sum_i a_ib_i$.

 This shows that $d(Q)$ is actually a homomorphism $I/I^2\xra{}A_*$.
 It is easy to check that the whole construction gives a homomorphism
 $d\:\Der(A)\xra{}\Hom_{A_*}(I/I^2,A_*)$.  If $d(Q)=0$ then all the
 elements $b_i$ are zero, so $Q\circ j_n=0$.  As $A$ is the homotopy
 colimit of the objects $A_n$, we conclude from
 Lemma~\ref{lem-lim-one} that $Q=0$.  Thus, $d$ is a monomorphism.
\end{proof}

The meaning of the proposition is elucidated by the following
elementary lemma.
\begin{lemma}\label{lem-I-sq}
 If $\{x_1,x_2,\ldots\}$ is a regular sequence in $R_*$, and $I$ is
 the ideal that it generates, then $I/I^2$ is freely generated over
 $R_*/I$ by the elements $x_i$.
\end{lemma}
\begin{proof}
 It is clear that $I/I^2$ is generated by the elements $x_i$.  Suppose
 that we have a relation $\sum_{i=1}^n a_ix_i=0$ in $I$ (not
 $I/I^2$).  We claim that $a_i\in I$ for all $i$.  Indeed, it is clear
 that $a_nx_n\in (x_1,\ldots,x_{n-1})$ so by regularity we have
 $a_n=\sum_{i=1}^{n-1}b_i x_i$ say; in particular, $a_n\in I$.
 Moreover, $\sum_{i=1}^{n-1}(a_i+b_ix_n)x_i=0$, so by induction we
 have $a_i+b_i x_n\in I$ for $i<n$, and thus $a_i\in I$ as required.

 Now suppose that we have a relation $\sum_i a_ix_i\in I^2$, say
 $\sum_i a_i x_i=\sum_{j\leq i}b_{ij}x_i x_j$.  We then have $\sum_i
 (a_i-\sum_{j\leq i}b_{ij}x_j)x_i=0$, so by the previous claim we have
 $a_i-\sum_{j\leq i}b_{ij}x_j\in I$, so $a_i\in I$.  This shows that
 the elements $x_i$ generate $I/I^2$ freely.
\end{proof}

\begin{corollary}\label{cor-self-map}
 In the situation of Proposition~\ref{prop-self-maps} the map
 $d\:\Der(A)\xra{}\Hom(I/I^2,A_*)$ is an isomorphism, and
 $\DD(A,A)_*$ is the completed exterior algebra generated by
 $\Der(A)$.  
\end{corollary}
\begin{proof}
 It is easy to see that $Q_i$ is a derivation and that
 $d(Q_i)(x_j)=\dl_{ij}$ (Kronecker's delta).  This shows that $d$ is
 surjective, and the rest follows. 
\end{proof}

\section{Formal group theory}\label{sec-formal}

In this section, we take $R=MU$, and let $F$ be the usual formal group
law over $MU_*$.  In places it will be convenient to use cohomological
gradings; we recall the convention $A^*=A_{-*}$.  We will write $q$
for the usual map $MU^*\xra{}BP^*$, and note that
$q(w_1)=v_1\pmod{2}$.

A well-known construction gives a power operation
\[ P\: R^dX \xra{} R^{2d}(\rp^2\tm X), \]
which is natural for spaces $X$ and strictly commutative ring spectra
$R$.  A good reference for such operations is~\cite{brmamcst:hir}; in
the case of $MU$, the earliest source is probably~\cite{td:sok}.  

In the case $R=S^{-1}MU$ there is an element $\ep\in R^2\rp^2$ such
that $R^*\rp^2=R^*[\ep]/(2\ep,\ep^2)$.  More generally, the
even-dimensional part of $R^*(\rp^2\tm X)$ is
$R^*(X)[\ep]/(2\ep,\ep^2)$, and $P(x)=x^2+\ep\tP(x)$ for a uniquely
determined operation $R^d(X)\xra{}R^{2d-2}(X)/2$.  We also have the
following properties:
\begin{align*}
 P(1)    &= 1                                   \\
 P(xy)   &= P(x) P(y)                           \\
 P(x+y)  &= P(x) + P(y) + (2+\ep w_1) xy        \\
 P(x)    &= x(x+_F\ep)
  \qquad \text{ if $x$ is the Euler class of a complex line bundle. }
\end{align*}

To handle the nonadditivity of $P$, we make the following
construction.  For any $MU^*$-algebra $A^*$, we define
\[ T(A^*) \eqdef 
    \{(r,s)\in A^*/2\tm A^*[\ep]/(2,\ep^2)\st s=r^2\pmod{\ep}\}.
\]
Given $a,b\in A^*$ (with $|b|=2|a|-2$) we define
$[a,b]\eqdef(a,a^2+\ep b)\in T(A^*)$.  We make $T(A^*)$ into a ring by
defining
\begin{align*}
 (r,s) + (t,u) &\eqdef (r+t, s+u+\ep w_1 r t) \\
 (r,s) . (t,u) &\eqdef (r t, s u)
\end{align*}
or equivalently
\begin{align*}
 [a,b] + [c,d] &\eqdef [a+c,b+d+w_1 a c] \\
 [a,b] . [c,d] &\eqdef [a c,a^2 d + b c^2]. 
\end{align*}
Note that $2[a,b]=[0,w_1 a^2]$ and $4[a,b]=0$, so $4T(A^*)=0$.  If we
define $Q(a)\eqdef(a,P(a))=[a,\tP(a)]$, then $Q$ gives a ring map
$MU^*X\xra{}T(MU^*X)$.

\begin{definition}\label{defn-induced}
 Suppose that $A^*$ is a PLRQ of $MU^*$, and let $f\:MU^*\xra{}A^*$ be
 the unit map.  We say that $A^*$ has an induced power operation (IPO)
 if there is a ring map $\Qb\:A^*\xra{}T(A^*)$ making the following
 diagram commute:
 \[ \diagram
   MU^* \dto_{f} \rto^{Q} & T(MU^*) \dto^{T(f)} \\
   A^* \rto_{\Qb} & T(A^*).
 \enddiagram \]
 Because $A^*$ is an LRQ, we know that such a map is unique if it
 exists.
\end{definition}

If $A^*=S^{-1}MU^*$ then we know that $S^{-1}MU$ can be constructed as a
strictly commutative $MU$-algebra and thus an $E_\infty$ ring spectrum,
and the power operation coming from this $E_\infty$ structure clearly
gives an IPO on $S^{-1}MU_*$.  For a more elementary proof, it suffices
to show that when $x\in S$ the image of $Q(x)$ in $T(S^{-1}MU_*)$ is
invertible.  However, the element $(x,x^2)$ is trivially invertible in
$T(S^{-1}MU)$ and $Q(x)$ differs from this by a nilpotent element, so
it too is invertible.  

It is now easy to reduce the following result to
Theorem~\ref{thm-even-general}. 
\begin{proposition}
 Let $A_*$ be a PLRQ of $MU_*$ which has no $2$-torsion and admits an
 IPO.  Then $A_*$ is strongly realisable.  \qed
\end{proposition}

We now give a formal group theoretic criterion for the existence of an
IPO.
\begin{definition}\label{defn-Zx}
 Let $F$ be a formal group law over a ring $A^*$.  Given an algebra
 $B^*$ over $A^*$ and an element $x\in B^2$, we define
 $Z(x)\eqdef(x,x(x+_F\ep))\in T(B^*)$.  (We need $x$ to be
 topologically nilpotent in a suitable sense to interpret this, but we
 leave the details to the reader.)  Thus, if $X$ is a space,
 $A^*=MU^*$, $B^*=MU^*X$ and $x$ is the Euler class of a complex line
 bundle over $X$ then $Z(x)=Q(x)$.
\end{definition}

\begin{proposition}\label{prop-IPO-FGL}
 Let $A^*$ be a LRQ of $MU^*$, and let $F$ be the obvious formal group
 law over $A^*$.  Then a ring map $\Qb\:A^*\xra{}T(A^*)$ is an IPO if
 and only if we have
 \[ Z(x) +_{\Qb_*F} Z(y) = Z(x+_Fy) \in T(\fps{A^*}{x,y}). \]
\end{proposition}
\begin{proof}
 Let $F'$ be the universal FGL over $MU^*$ and put
 $Z'(x)=(x,x(x+_{F'}\ep))$.  Let $f\:MU^*\xra{}A^*$ be the unit map,
 so that $F=f_*F'$.  Using the universality of $F'$, we see that $\Qb$
 is an IPO if and only if we have
 \[ X+_{\Qb_*f_*F'}Y = X +_{T(f)_*Q_*F'}Y \in \fps{T(A^*)}{X,Y}. \] 
 The left hand side is of course $X+_{\Qb_*F}Y$. There is an evident
 map
 \[ \fps{T(A^*)}{X,Y} \xra{} T(\fps{A^*}{x,y}), \]
 sending $X$ to $Z(x)$ and $Y$ to $Z(y)$, and one can check that this
 is injective.  Thus, $\Qb$ is an IPO if and only if
 \[ Z(x) +_{\Qb_*F} Z(y) = Z(x)+_{T(f)_*Q_*F'}Z(y)
     \in T(\fps{A^*}{x,y}). 
 \]
 The right hand side here is $T(f)(Z'(x)+_{Q_*F'}Z'(y))$ and
 $Z(x+_Fy)=T(f)(Z'(x+_{F'}y))$ so the proposition will follow once we
 prove that $Z'(x+_{F'}y)=Z'(x)+_{Q_*F'}Z'(y)\in T(\fps{MU^*}{x,y})$.
 To do this, we use the usual isomorphism
 $\fps{MU^*}{x,y}=MU^*(\cpi\tm\cpi)$, so that $x$, $y$ and $x+_{F'}y$
 are Euler classes, so $Q(x)=Z'(x)$ and $Q(y)=Z'(y)$ and
 $Q(x+_{F'}y)=Z'(x+_{F'}y)$.  As $Q$ is a natural multiplicative
 operation we also have
 $Q(x+_{F'}y)=Q(x)+_{Q_*F'}Q(y)=Z(x)+_{Q_*F'}Z(y)$, which gives the
 desired equation.
\end{proof}

We now use this to show that there is an IPO on $kU^*$.  In this case
the real reason for the IPO is that the Todd genus gives an $H_\infty$
map $MU\xra{}kU$, but we give an independent proof as a warm-up for the
case of $BP^*$.

\begin{proposition}\label{prop-IPO-kU}
 Let $f\:MU^*\xra{}kU^*\eqdef\Zh[u]$ be the Todd genus.  Then there is
 an induced power operation on $kU^*$, given by $\Qb(u)=[u,u^3]$.
 Thus, $kU^*$ and $KU^*$ are strongly realisable.
\end{proposition}

\begin{proof}
 The FGL over $kU^*$ coming from $f$ is just the multiplicative FGL
 $x+_Fy=x+y+uxy$, so $Z(x)=(x,x^2+x\ep+ux^2\ep)=[x,x+ux^2]$.  If we
 put $U=[u,u^3]$ then $X+_{\Qb_*F}Y=X+Y+UXY$.  We thus need only
 verify that $Z(x)+Z(y)+U\,Z(x)Z(y)=Z(x+y+uxy)$.  This is a
 straightforward calculation; some steps are as follows.
 \begin{align*}
  Z(x) + Z(y) &= [x+y,x+y+u(x^2+xy+y^2)]                \\
  Z(x)Z(y)    &= [xy,xy(x+y)]                           \\
  U\,Z(x)Z(y) &= [uxy,u^2xy(x+y)+u^3x^2y^2]             \\
  Z(x+y+uxy)  &= [x+y+uxy,x+y+u(x^2+xy+y^2)+u^3x^2y^2].
 \end{align*}
\end{proof}

We now turn to the case of $BP^*$.  For the moment we prove only that
an IPO exists; in the next section we will calculate it.

\begin{proposition}\label{prop-IPO-BP}
 There is an IPO on $BP^*$, so $BP^*$ is strongly realisable.
\end{proposition}

This is proved after Lemma~\ref{lem-zxqf}.

\begin{definition}\label{defn-z}
 For the rest of this section, we will write
 \[ z \eqdef z(x) \eqdef 
    \sum_{k\ge 0} v_1^{2^k}x^{2^k} \in \fps{\FF_2[v_1]}{x}.
 \]
 Note that 
 \[ z^2 = z + v_1 x \]
 so
 \[ z/(v_1 x) = 1/(1 + z). \]
\end{definition}

\begin{lemma}\label{lem-invdif}
 We have
 \[ x +_F \ep = x + (1 + z)\ep  \qquad \text{ in } \qquad
    \fps{BP^*[\ep]}{x}/(2\ep,\ep^2).
 \]
\end{lemma}
\begin{proof}
 Working rationally and modulo $\ep^2$, we have
 \[ \log_F(\ep) = \ep \]
 so
 \[ x +_F\ep = \exp_F(\log_F(x)+\ep) = x + \exp'_F(\log_F(x))\ep =
    x + \ep/\log'_F(x).
 \]
 Note that $\log'_F(x)$ is integral and its constant term is $1$, so
 the above equation is between integral terms and we can sensibly
 reduce it modulo $2$.

 We next recall the formula for $\log_F(x)$ given
 in~\cite[Section~4.3]{ra:ccs}.  We consider sequences
 $I=(i_1,\ldots,i_r)$ with $l\geq 0$ and $i_j>0$ for each $j$.  We
 write $|I|\eqdef r$ and $\|I\|\eqdef i_1+\cdots+i_r$.  We also write
 \[ v_I \eqdef v_{i_1}^{m_1}\ldots v_{i_r}^{m_r} \]
 where
 \[ m_j \eqdef 2^{\sum_{k<j} i_k}. \]
 The formula is
 \[ \log_F(x) = \sum_I v_I x^{2^{\|I\|}}/2^{|I|}. \]
 The only terms which contribute to $\log'_F(x)$ modulo $2$ are those
 for which $\|I\|=|I|$, so $i_j=1$ for all $j$.  If $I$ has this form
 and $|I|=k$ then $v_I=v_1^{2^k-1}$.  Thus
 \[ \log'_F(x) = \sum_k v_1^{2^k-1} x^{2^k-1} = z/(v_1 x) \pmod{2}. \]
 As remarked in Definition~\ref{defn-z}, we have $z/(v_1 x)=1/(1+z)$,
 so 
 \[ x+_F\ep = x + \ep/\log'_F(x) = x + (1+z)\ep \]
 as claimed.
\end{proof}

\begin{lemma}\label{lem-zxqf}
 In $T(\fps{BP^*}{x,y})$ we have
 \[ [0,x] +_{QF} [0,y] = [0,x+y] \]
 and
 \[ [x,y] = [x,0] +_{QF} 
      \lb 0,\sum_{k\ge 0} (v_1 x)^{2(2^k-1)} y\rb.
 \]
 In particular, we have
 \[ Z(x) = [x,0] +_{QF} [0,z/v_1] =
      [x,0] +_{QF} \lb 0,\sum_{k\ge 0} v_1^{2^k-1}x^{2^k}\rb.
 \]
\end{lemma}
\begin{proof}
 The first statement is clear, just because $[0,x][0,y]=0$.  For the
 second statement, write $X\eqdef[x,0]$ and 
 \[ w\eqdef \sum_k(v_1x)^{2(2^k-1)}y = y(z/v_1 x)^2 = y/(1+z)^2, \]
 and $W\eqdef [0,w]$.  Let $a_{ij}\in MU_{2(i+j-1)}$ be the
 coefficient of $x^iy^j$ in $x+_Fy$.  Because $W^2=0$ we have
 $X+_{QF}W=X+W+\sum_{j>0}Q(a_{1j})X^jW$, and
 \[ Q(a_{1j})X^j W =
     [a_{1j},\tP(a_{1j})][x^j,0][0,w] = [0,a^2_{1j}x^{2j}w].
 \]
 This expression is to be interpreted in $T(\fps{BP^*}{x,y})$, so we
 need to interpret $a_{1j}$ in $BP^*/2$.  Thus Lemma~\ref{lem-invdif}
 tells us that $a_{10}=1$ and $a_{1,2^k}=v_1^{2^k}$ and all other
 $a_{1j}$'s are zero.  Thus
 \begin{multline*}
   X+_{QF}W = X + W + \sum_{k\ge 0} [0,(v_1 x)^{2^{k+1}}w] = \\
     \lb x,w\left(1+\sum_{k\ge 0} (v_1 x)^{2^{k+1}}\right)\rb
      = [x,w(1+z^2)] = [x,y]
 \end{multline*}
 as claimed.  

 For the last statement, Lemma~\ref{lem-invdif} gives 
 \[ Z(x) = (x,x(x+_F\ep)) = (x,x^2 + x(1+z)\ep) = [x,x(1+z)]. \]
 By the previous paragraph, this can be written as
 $[x,0]+_{QF}[0,x(1+z)/(1+z)^2]=[x,0]+_{QF}[0,z/v_1]$.
\end{proof}

\begin{proof}[Proof of Proposition~\ref{prop-IPO-BP}]
 To show that $\Qb$ exists, it is enough to show that the formal group
 law on $T(BP^*)$ obtained from the map
 $MU^*\xra{Q}T(MU^*)\xra{T(q)}T(BP^*)$ is $2$-typical.  Let $p$ be an
 odd prime, so the associated cyclotomic polynomial is
 $\Phi_p(t)=1+t+\cdots+t^{p-1}$.  We need to show that
 \[ X +_{QF} \Om X +_{QF} \cdots +_{QF} \Om^{p-1} X = 0 
    \qquad\text{ in }\qquad B^*\eqdef\fps{T(BP^*)[\Om]}{X}/\Phi_p(\Om).
 \]
 (This is just the definition of $2$-typicality for formal groups over
 rings which may have torsion.)  Consider the ring
 $C^*\eqdef T(\fps{BP^*[\om]}{x}/\Phi_p(\om))$.  As $\Phi_p(\om)=0$ we
 have $t^p-1=\prod_{j=0}^{p-1}(t-\om^j)$, and by looking at the
 coefficient of $t^{p-2}$ we find that $\sum_{0\le i<j<p}\om^{i+j}=0$.
 Now write $\Om\eqdef[\om,0]$ and $X\eqdef[x,0]$, so that
 $\Om,X\in C^*$.  We find that 
 \[ \Phi_p(\Om) = \sum_{i=0}^{p-1} [\om^i,0] =
    \lb \Phi_p(\om),v_1\sum_{0\le i<j<p}\om^{i+j}\rb = 0.
 \]
 This gives us a ring map $B^*\xra{}C^*$; we claim that this is
 injective.  Indeed, it is easy to see that
 $\{\om,\om^2,\ldots,\om^{p-1}\}$ is a basis for
 $\Zh[\om]/\Phi_p(\om)$, and that $\al\mapsto\al^2$ is a 
 permutation of this basis.  Suppose that we have
 \[ \sum_{i=1}^{p-1}\sum_{j\ge 0} [a_{ij},b_{ij}]\Om^iX^j=0
    \text{ in } C^*.
 \]
 Using the evident map $C^*\xra{}\fps{BP^*[\om]}{x}/(2,\Phi_p(\om))$,
 we see that $a_{ij}=0$ for all $i,j$.  As
 $[0,b]\Om^iX^j=[0,\om^{2i}X^{2j}b]$, we see that
 \[ \sum_{i=1}^{p-1}\sum_{j\ge 0} b_{ij}\om^{2i}X^{2j}=0.
 \]
 As the elements $\om^{2i}$ are a permutation of the elements $\om^i$,
 we see that $b_{ij}=0$ for all $i,j$.  We may thus regard $B^*$ as a
 subring of $C^*$.

 Next, we know that
 \begin{equation}\label{eqn-ind-BP-A}
   Z(x) +_{QF} Z(\om x) +_{QF}\cdots+_{QF} Z(\om^{p-1}x) = 
    Z(x+_F \om x +_F \cdots +_F \om^{p-1}x) = 0,
 \end{equation}
 because $F$ is $2$-typical over $BP_*$.  By Lemma~\ref{lem-zxqf}, we
 also know that 
 \begin{equation}\label{eqn-ind-BP-B}
  Z(\om^ix)=\Om^iX+_{QF}[0,w_i], 
 \end{equation}
 where $w_i = \sum_k v_1^{2^k-1} \om^{2^ki} x^{2^k}$.
 It is easy to see that $[0,w_i][0,w_j]=0$, so that
 $[0,w_i]+_{QF}[0,w_j]=[0,w_i+w_j]$.  We also have
 $\sum_{i=0}^{p-1}\om^{2^ki}=0$ for all $k$.  This means that
 \begin{equation}\label{eqn-ind-BP-C}
   \sum^{QF}_i [0,w_i] = 
    [0,\sum_k v_1^{2^k-1}x^{2^k} \sum_{i=0}^{p-1}\om^{2^ki}]=0.
 \end{equation}
 By combining equations~(\ref{eqn-ind-BP-A}) to~(\ref{eqn-ind-BP-C}),
 we see that
 \[ \sum^{QF}_i \Om^i X = 0 \]
 as required.
\end{proof}

\section{The power operation on $BP^*$}\label{sec-Q-BP}

We now give explicit formulae for the IPO on $BP^*$.

\begin{definition}\label{defn-u}
 Given a subset $J=\{j_1<\ldots<j_r=n\}\sse\{1,\ldots,n\}$, we define
 \[ u_J \eqdef
     v_{j_1+1}\prod_{k=1}^{r-1}
       (v_1 v_{j_k})^{2(2^{j_{k+1}-j_k}\,-1)}
     \in\pi_{|v_{n+1}|}BP
 \]
 and $u_n=\sum_Ju_J$, where $J$ runs over subsets of $\{1,\ldots,n\}$
 that contain $n$.  By separating out the case $r=1$ and putting
 $j=n-j_{r-1}$ in the remaining cases we obtain a recurrence relation
 \[ u_n = v_{n+1} + \sum_{j=1}^{n-1} (v_1v_{n-j})^{2(2^j-1)}u_{n-j}.
 \]
\end{definition}

\begin{proposition}\label{prop-Q-formula}
 The induced power operation on $BP^*$ is given by 
 \[ \Qb(v_n) = \begin{cases}
      {}[0,v_1]                 & \text{ if } n=0               \\
      {}[v_1,v_2]               & \text{ if } n=1               \\
      {}[v_n,v_1 v_n^2 + u_n]   & \text{ if } n>1
    \end{cases}
 \]
 Moreover, we have $u_n=v_{n+1}\pmod{v_1^2}$.
\end{proposition}
This is proved after Corollary~\ref{cor-exp-QF}.  We will reuse the
notation of Definition~\ref{defn-z}.

\begin{lemma}\label{lem-exp-F}
 We have $\exp_F(2x)=2z/v_1$ in $\fps{BP^*}{x}/4$.
\end{lemma}
\begin{proof}
 Using Ravenel's formulae as in the proof of Lemma~\ref{lem-invdif},
 we have
 \[ \log_F(2x)/2 = \sum_I 2^{2^{\|I\|}-|I|-1}v_Ix^{2^{\|I\|}}. \]
 When $k\ge 0$ we have $2^k\ge k+1$, with equality only when $k=0$ or
 $k=1$.  It follows easily that 
 \[ \log_F(2x)/2 = x + v_1 x^2 \pmod{2}. \]
 By inverting this, we find that
 \[ \exp_F(2x)/2 = \sum_{k\ge 0} v_1^{2^k-1}x^{2^k}=z/v_1 \pmod{2}, \]
 and thus that $\exp_F(2x)=2z/v_1\pmod{4}$.
\end{proof}

Because $T(BP^*)$ is a torsion ring, the formal group law $QF$ has no
$\exp$ series.  Nonetheless, $\exp_F(2X)$ is a power series over
$BP^*$, so we can apply $Q$ to the coefficients to get a power series
over $T(BP^*)$ which we call $\exp_{QF}(2X)$.  This makes perfect
sense even though $\exp_{QF}(X)$ does not.
\begin{corollary}\label{cor-exp-QF}
 In $\fps{T(BP^*)}{X}$, we have
 \[ \exp_{QF}(2X) = \sum_{k\ge 0} [0,v_1^{2^{k+1}-1}] X^{2^k}. \]
 By taking $X=Z(x)\in T(\fps{BP^*}{x})$, we get
 \[ \exp_{QF}(2Z(x)) = \lb 0,\sum_{j>0} v_1^{2^j-1} x^{2^j}\rb
                   = [0, z/v_1 + x] .
 \]
\end{corollary}
\begin{proof}
 Because $4=0$ in $T(BP^*)$, it follows immediately from the lemma
 that $\exp_{QF}(2X)=\sum_k 2\Qb(v_1)^{2^k-1}X^{2^k}$.  Using
 $\Qb(v_1)=[v_1,\tP(v_1)]$, we see that
 $2\Qb(v_1)^{2^k-1}=[0,v_1^{2^{k+1}-1}]$, and the first claim
 follows.  If we now put $X=Z(x)=[x,x(1+z)]$ then
 $[0,v_1^{2^{k+1}-1}]X^{2^k}=[0,v_1^{2^{k+1}-1}x^{2^{k+1}}]$, and the
 second claim follows.
\end{proof}

\begin{proof}[Proof of Proposition~\ref{prop-Q-formula}]
 Let $p_k$ denote the image of $\tP(v_k)\in MU^{2^{k+2}-2}/2$ in
 $BP^*/2$ and write $V_k=\Qb(v_k)=[v_k,p_k]\in T(BP^*)$.  Recall that
 the Hazewinkel generators $v_k$ are characterised by the formula
 \[ [2]_F(x)=\exp_F(2x)+_F\sum^F_{k>0}v_kx^{2^k} \in\fps{BP^*}{x}. \]
 By applying the ring map $\Qb$ and putting $X=Z(x)$ we obtain
 \[ [2]_{QF}(Z(x)) = 
    \exp_{QF}(2Z(x)) +_{QF} \sum^{QF}_{k>0} V_k Z(x)^{2^k}
    \in T(\fps{BP^*}{x}).
 \]
 The first term can be evaluated using Corollary~\ref{cor-exp-QF}.
 For the remaining terms, we have 
 \[ V_k Z(x)^{2^k} = [v_k,p_k][x^{2^k},0] =
     [v_k x^{2^k},p_k x^{2^{k+1}}].
 \]
 We can use Lemma~\ref{lem-zxqf} to rewrite this as
 \begin{align*}
   V_k Z(x)^{2^k} & = 
     [v_k x^{2^k},0] +_{QF} 
      \lb 0,\sum_{l>0}(v_1v_k x^{2^k})^{2^l-2}p_k x^{2^{k+1}}\rb \\
   &=
     [v_k x^{2^k},0] +_{QF} 
      \lb 0,\sum_{l>0} (v_1v_k)^{2^l-2} p_k x^{2^{k+l}}\rb.
 \end{align*}
 After using the formula $[0,b]+_{QF}[0,c]=[0,b+c]$ to collect terms,
 we find that
 \begin{equation}\label{eqn-Q-formula-A}
   [2]_{QF}(Z(x)) = 
     \lb 0,\sum_{l>0} v_1^{2^l-1}x^{2^l} + 
        \sum_{k,l>0} (v_1v_k)^{2^l-2} p_k x^{2^{k+l}}\rb +_{QF}
     \sum^{QF}_{k>0} [v_k x^{2^k},0].
 \end{equation}
 On the other hand, we know that 
 \begin{align*}
  [2]_{QF}(Z(x)) &= Z([2]_F(x))                                 \\
                 &= Z\left(\exp_F(2x)+_F\sum^F_{k>0}v_kx^{2^k}\right)\\
                 &= Z(\exp_F(2x))+_{QF}\sum^{QF}_{k>0}Z(v_kx^{2^k}).
 \end{align*}
 The first term is zero because $\exp_F(2x)$ is divisible by $2$.  For
 the remaining terms, Lemma~\ref{lem-zxqf} gives
 \[ Z(v_kx^{2^k}) =
     [v_k x^{2^k},0] +_{QF} 
      [0,\sum_{j\ge 0} v_1^{2^j-1} v_k^{2^j} x^{2^{k+j}}].
 \]
 Thus, we have
 \begin{equation}\label{eqn-Q-formula-B}
   [2]_{QF}(Z(x)) = 
     \lb 0,\sum_{k>0}\sum_{l\ge 0} v_1^{2^l-1} v_k^{2^l} x^{2^{k+l}}\rb
     +_{QF} \sum^{QF}_{k>0} [v_k x^{2^k},0].
 \end{equation}
 By comparing this with equation~(\ref{eqn-Q-formula-A}) and equating
 coefficients of $x^{2^{n+1}}$, we find that
 \[ v_1^{2^{n+1}-1} + \sum_{j=1}^n (v_1v_{n+1-j})^{2^j-2} p_{n+1-j} =
    \sum_{j=0}^n v_1^{2^j-1}v_{n+1-j}^{2^j}.
 \]
 After some rearrangement and reindexing, this becomes
 \[ p_n + v_1 v_n^2 = v_1^{2^{n+1}-1} + v_{n+1} +
     \sum_{j=1}^{n-1} (v_1v_{n-j})^{2(2^j-1)}(p_{n-j}+v_1v_{n-j}^2).
 \]
 In particular, we have $p_1=v_2$.  We now define
 \[ p'_n = \begin{cases}
      v_1                       & \text{ if } n=0               \\
      v_2                       & \text{ if } n=1               \\
      v_1 v_n^2 + u_n           & \text{ if } n>1
    \end{cases}
 \]
 The claim of the proposition is just that $p_n=p'_n$ for all
 $n\ge 0$.  Using the recurrence relation given in
 definition~\ref{defn-u}, one can check that for all $n>0$ we have
 \[ p'_n + v_1 v_n^2 = v_1^{2^{n+1}-1} + v_{n+1} +
     \sum_{j=1}^{n-1} (v_1v_{n-j})^{2(2^j-1)}(p'_{n-j}+v_1v_{n-j}^2).
 \]
 In particular, we have $p'_1=v_2=p_1$, and it follows inductively
 that $p_n=p'_n$ for all $n>0$.  We also have 
 \[ Q(v_0)=Q(1)+Q(1)=[1,0]+[1,0]=[0,v_1] \]
 so $p_0=v_1=p'_0$.
\end{proof}

\begin{remark}\label{rem-Pvk}
 The first few cases are
 \begin{align*}
  p_0 &= v_1                                               \\
  p_1 &= v_2                                               \\
  p_2 &= v_1^4 v_2 + v_1 v_2^2 + v_3                       \\
  p_3 &= v_1^{12} v_2 + v_1^6 v_2^3 + v_1^2 v_2^2 v_3 +
          v_1 v_3^2 + v_4.
 \end{align*}
 In particular, we find that $p_3\not\in(v_k\st k\geq 3)$, which shows
 that there is no commutative product on $BP\langle 2\rangle$,
 considered as an object of $\DD$.  This problem does not go away if
 we replace the Hazewinkel generator $v_k$ by the corresponding Araki
 generator, or the bordism class $w_k$ of a smooth quadric
 hypersurface in $\mathbb{C}P^{2^k}$.  However, it is possible to
 choose a more exotic sequence of generators for which the problem
 does go away, as indicated by the next result.
\end{remark}

\begin{proposition}\label{prop-ideal-J}
 Fix an integer $n>0$.  There is an ideal $J\leq BP^*$ such that
 \begin{enumerate}
  \item The evident map
   \[ \Ztl[v_1,\ldots,v_n] \xra{} BP^* \xra{} BP^*/J \]
   is an isomorphism.
  \item $\tP(J)\leq J\pmod{2}$.
  \item $I_n+J=I_n+(v_k\st k>n)=(v_k\st k\neq n)$.
 \end{enumerate}
 The proof will construct an ideal explicitly, but it is not the only
 one with the stated properties.  If $n=1$ we can take
 $J=(v_k\st k>n)$, but for $n>1$ this violates condition~(2).
\end{proposition}
\begin{remark}
 The subring $\Ztl[v_1,\ldots,v_n]$ of $BP^*$ is the same as the
 subring generated by all elements of degree at most $2^{n+1}-2$; it
 is thus defined independently of the choice of generators for
 $BP_*$. 
\end{remark}
\begin{proof}
 First consider the case $n=1$, and write $J=(v_k\st k>1)$.  By
 inspecting definition~\ref{defn-u}, we see that $u_n\in J$ for all
 $n>1$, and thus Proposition~\ref{prop-Q-formula} tells us that
 $\tP(J)\leq J\pmod{2}$.  We may thus assume that $n>1$.  Write
 $B^*=\Ztl[v_1,\ldots,v_n]$, thought of as a subring of $BP^*$.  We
 will recursively define a sequence of elements $x_k\in BP^*$ for
 $k>n$ such that
 \begin{itemize}
  \item[(a)] $x_k\in v_k + v_1^2 B^*$
  \item[(b)] $\tP(x_{k-1})\in(x_{n+1},\ldots,x_k)\pmod{2}$ \;if\;
   $k>n+1$.
 \end{itemize}
 It is clear that we can then take $J=(x_k\st k>n)$.  We start by
 putting $x_{n+1}=v_{n+1}$.  Suppose that we have defined
 $x_{n+1},\ldots,x_k$ with the stated properties.  There is an evident
 map
 \[ \FF[v_1,\ldots,v_n,v_{k+1}] \xra{f} BP_*/(2,x_{n+1},\ldots,x_k),\]
 which is an isomorphism in degree $2(2^{k+1}-1)=|v_{k+1}|$.  Let
 $\pb_k$ be the image of $\tP(x_k)$ in $BP_*/(2,x_{n+1},\ldots,x_k)$,
 and write $\xb_{k+1}=f^{-1}(\pb_k)$.  We can lift this to get an
 element $x_{k+1}$ of $\Ztl[v_1,\ldots,v_n,v_{k+1}]$ such that
 $\xb_{k+1}=x_{k+1}\pmod{2}$ and every coefficient in $x_{k+1}$ is $0$
 or $1$.  It is easy to see that condition~(b) is satisfied, and that
 $x_{k+1}\in v_{k+1}+B^*$.  However, we still need to show that
 $x_{k+1}-v_{k+1}$ is divisible by $v_1^2$.  By assumption we have
 $x_k=v_k+v_1^2b$ for some $b\in B^*$.  Recall from
 Proposition~\ref{prop-Q-formula} that
 $\tP(v_k)=v_{k+1}+v_1v_k^2\pmod{2,v_1^2}$.  It follows after a small
 calculation that $\tP(x_k)=v_{k+1}+v_1v_k^2\pmod{2,v_1^2}$ also.
 Moreover, we have $v_k^2=v_1^4 b^2\pmod{2,x_k}$, so
 $\pb_k=v_{k+1}\pmod{2,v_1^2}$.  It follows easily that
 $x_{k+1}=v_{k+1}\pmod{v_1^2}$, as required.
\end{proof}

We give one further calculation, closely related to
Proposition~\ref{prop-Q-formula}.  
\begin{proposition}\label{prop-Pwk}
 Recall that $I_k\eqdef(w_0,\ldots,w_{k-1})<MU^*$, where $w_i$ is the
 bordism class of a smooth quadric hypersurface in $\cp^{2^i}$.  We
 have $\tP(I_{k-1})\leq I_k$, and $\tP(w_{k-1})=w_k\pmod{I_k}$.
\end{proposition}
\begin{proof}
 If $k=1$ we have $I_0=0$ and $w_0=2$, so
 $P(w_0)=P(1)+P(1)+(2+w_1\ep)=w_1\ep\pmod{2}$, as required.  Thus, we
 may assume that $k>1$, and it follows easily from the formulae for
 $P(x+y)$ and $P(xy)$ that $P$ induces a ring map
 $MU^*\xra{}B^*=(MU^*/I_k)[\ep]/\ep^2$.  Note that
 $[2]_F(x)=w_kx^{2^k}+O(x^{2^k+1})$ over $B^*$.  Write
 $X=x(x+_F\ep)\in\fps{MU^*[\ep]}{x}/(I_k,\ep^2)$.  Arguing in the
 usual way, we see that 
 \[ [2]_{P_*F}(X)= [2]_F(x)([2]_F(x)+_F\ep) =
     \ep w_k x^{2^k} + O(x^{2^k+1}). 
 \]
 It follows easily that we must have
 \[ [2]_{P_*F}(X) = \ep w_k X^{2^{k-1}} + O(X^{2^{k-1}+1}). \]
 It follows that $P(w_i)=0\in B^*$ for $i<k-1$, and that
 $P(w_{k-1})=\ep w_k\in B^*$, as required.
\end{proof}

\section{Applications to $MU$}\label{sec-mu}

\begin{proof}[Proof of Proposition~\ref{prop-MU-omni}]
 The claims involving $kU$ and $KU$ follow from
 Proposition~\ref{prop-IPO-kU}, and those for $BP$ follow from
 Proposition~\ref{prop-IPO-BP}.  The claim $H$ follows from
 Theorem~\ref{thm-even-general}, as the condition
 $\tP(I)\sse I\pmod{2}$ is trivially satisfied for dimensional
 reasons.  The claim for $H\FF$ can be proved in the same way as
 Theorem~\ref{thm-even-proof} after noting that all the obstruction
 groups are trivial.
\end{proof}

\begin{proof}[Proof of Proposition~\ref{prop-BPn}]
 Choose an ideal $J$ as in Proposition~\ref{prop-ideal-J} and set
 $\BP{n}'_*=BP_*/J$.  Everything then follows from
 Theorem~\ref{thm-even-general}.  
\end{proof}

We now take $R=MU_{(2)}$ and turn to the proof of
Theorem~\ref{thm-Pn}.  As previously, we let $w_k\in\pi_{2^{k+1}-2}R$
denote the bordism class of the quadric hypersurface $W_{2^k}$ in
$\cp^{2^k}$.  Recall that the image of $w_k$ in $BP_*$ is $v_k$ modulo
$I_k=(v_0,\ldots,v_{k-1})$, and thus
$P(n)_*=BP_*/(w_0,\ldots,w_{n-1})$.

We next choose a product $\phi_k$ on $MU/w_k$ for each $k$.  For $k=0$
we choose there are two possible products, and we choose one of them
randomly.  (It is possible to specify one of them precisely using
Baas-Sullivan theory, but that would lead us too far afield.)  For
$k>0$, we recall from Proposition~\ref{prop-Pwk} that
$\tP(w_k)=w_{k+1}\pmod{I_{k+1}}$.  It follows easily that there is a
product $\phi_k$ such that $c(\phi_k)=w_{k+1}\pmod{w_1,\ldots,w_k}$,
and that this is unique up to a term $u\circ(\btt)$ with
$u\in(w_1,\ldots,w_k)$.  From now on, we take $\phi_k$ to be a product
with this property.  It is easy to see that the resulting product
$MU/w_0\Smash\ldots\Smash MU/w_{n-1}$ is independent of the choice of
$\phi_k$'s (except for $\phi_0$).

\begin{definition}\label{defn-MUIn}
 We write
 \[ MU/I_n = MU/w_0\Smash\ldots\Smash MU/w_{n-1}, \]
 made into a ring as discussed above.  For $i<n$, we define
 \[ Q_i\: MU/I_n \xra{} \Sg^{2^{i+1}-1} MU/I_n \]
 by smashing the Bockstein map $\btb\:MU/w_i\xra{}\Sg^{2^{i+1}-1}MU/w_i$
 with the identity on the other factors.  We also define
 \begin{align*}
  P(n) &\eqdef BP\Smash MU/I_n                \\
  B(n) &\eqdef w_n^{-1}P(n)                  \\
  k(n) &\eqdef \BP{n}'\Smash MU/I_n           \\
  K(n) &\eqdef w_n^{-1}k(n).
 \end{align*}
 It is clear that $\pi_*(MU/I_n)=MU_*/I_n$ and $\pi_*(P(n))=P(n)_*$ and
 $\pi_*(B(n))=B(n)_*$.  Condition~(2) in Proposition~\ref{prop-BPn}
 assures us that $\pi_*k(n)=k(n)_*$ and $\pi_*K(n)=K(n)_*$ as well.
 As $BP$ and $\BP{n}'$ are commutative, it is easy to see that $P(n)$,
 $B(n)$, $k(n)$ and $K(n)$ are central algebras over $BP$,
 $v_n^{-1}BP$, $\BP{n}'$ and $E(n)'$ respectively.  The derivations
 $Q_i$ on $MU/I_n$ clearly induce compatible derivations on $P(n)$,
 $B(n)$, $k(n)$ and $K(n)$.
\end{definition}

\begin{proposition}\label{prop-twist-Q}
 The product $\phi$ on $MU/I_n$ satisfies
 \[ \phi-\phi\circ\tau = w_n \phi\circ(Q_{n-1}\Smash Q_{n-1}). \]
 Similarly for $P(n)$, $B(n)$, $k(n)$ and $K(n)$.
\end{proposition}
\begin{proof}
 This follows easily from the fact that $c(\phi_{k-1})=w_k\pmod{I_k}$,
 given by Proposition~\ref{prop-Pwk}.
\end{proof}

\begin{proposition}\label{prop-Pn-map}
 Let $A$ be a central $BP$-algebra such that $\pi_0(A)=\{0,1\}$,
 $\pi_{2^{n+1}-2}(A)=\{0,v_n\}$ and $\pi_k(A)=0$ for
 $0<k<2^{n+1}-2$.  Then either there is a unique map $P(n)\xra{}A$ of
 $BP$-algebras, or there is a unique map $P(n)\xra{}A^{\text{op}}$
 (but not both).  Analogous statements hold for $B(n)$, $k(n)$ and
 $K(n)$ with $BP$ replaced by $v_n^{-1}BP$, $\BP{n}'$ and $E(n)'$
 respectively.  
\end{proposition}
\begin{proof}
 We treat only the case of $P(n)$; the other cases are essentially
 identical.  Any ring map $MU/I_n\xra{}A$ commutes with the given map
 $BP\xra{}A$, because the latter is central.  It follows that maps
 $P(n)\xra{}A$ of $BP$-algebras biject with maps $MU/I_n\xra{}A$ of
 rings, which biject with systems of commuting ring maps
 $MU/w_i\xra{}A$ for $0\leq i<n$.  For $i<n-1$ we have
 $\pi_{2|w_i|+2}(A)=0$, so Proposition~\ref{prop-maps-Rx} tells us
 that the unique unital map $f_i\:MU/w_i\xra{}A$ is a ring map.  This
 remains the case if we replace the product $\psi$ on $A$ by
 $\psi\circ\tau$, or in other words replace $A$ by $A^{\text{op}}$.
 There is an obstruction
 $d_A(\phi_{n-1})\in\pi_{2^{n+1}-2}(A)=\{0,v_n\}$ which may prevent
 $f_{n-1}$ from being a ring map.  If it is nonzero, we have
 \[ d_A(\phi_{n-1}\circ\tau) = 
    d_A(\phi_{n-1}+\tP(w_{n-1})\circ(\btt)) = 
    d_A(\phi_{n-1}) + v_n = 0
 \]
 This shows that $f_{n-1}\:MU/w_{n-1}\xra{}A^{\text{op}}$ is a ring
 homomorphism.  After replacing $A$ by $A^{\text{op}}$ if necessary,
 we may thus assume that all the $f_i\:MU/w_i\xra{}A$ are ring maps.

 The obstruction to $f_i$ commuting with $f_j$ lies in
 $\pi_{|w_i|+|w_j|+2}(A)$.  If $i$ and $j$ are different then at least
 one is strictly less than $n-1$; it follows that
 $|w_i|+|w_j|+2<2^{n+1}-2$ and thus that the obstruction group is
 zero.  Thus $f_i$ commutes with $f_j$ when $i\neq j$, and we get a
 unique induced map $MU/I_n\xra{}A$, as required.
\end{proof}

\section{Point-set level foundations}\label{sec-found}

In order to analyse the commutativity obstruction $\cb(x)$ more
closely and relate them to power operations, we need to recall some
internal details of the EKMM category.

EKMM use the word ``spectrum'' in the sense defined by Lewis and
May~\cite{lemast:esh}, rather than the sense we use elsewhere in this
paper.  They construct a category $\LL\CS$ of ``$\LL$-spectra''.  This
depends on a universe $\CU$, but the functor
$\CL(\CU,\CV)\ths_{\CL(U)}(-)$ gives a canonical equivalence of
categories from $\LL$-spectra over $\CU$ to $\LL$-spectra over $\CV$,
so the dependence is only superficial.  (Here $\CL(\CU,\CV)$ is the
space of linear isometries from $\CU$ to $\CV$.)  We therefore take
$\CU=\rinf$.  EKMM show that $\LL\CS$ has a commutative and
associative smash product $\Smash_\CL$, which is not unital.  However,
there is a sort of ``pre-unit'' object $S$, with a natural map
$S\Smash_\CL X\xra{}X$.  They then define the subcategory
$\CM\eqdef\CM_S=\{X\st S\Smash_\CL X=X\}$ of ``$S$-modules'', and
prove that $S\Smash_\CL S=S$ so that $S\Smash_\CL X$ is an $S$-module
for any $X$.  We write $\Smash$ for the restriction of $\Smash_\CL$ to
$\CM$.

We next give a brief outline of the properties of $\CM$.  Let $\CT$ be
the category of based spaces (all spaces are assumed to be compactly
generated and weakly Hausdorff).  We write $0$ for the one-point
space, or for the basepoint in any based space, or for the trivial map
between based spaces.

We give $\CT$ the usual Quillen model structure for which the
fibrations are Serre fibrations.  We write $\hT$ for the category with
Hom sets $\pi_0F(A,B)=\CT(A,B)/\text{homotopy}$, and $\hbT$ for the
category obtained by inverting the weak equivalences.  We refer to
$\hT$ as the strong homotopy category of $\CT$, and $\hbT$ as the weak
homotopy category.

The category $\CM$ is a topological category: the Hom sets $\CM(X,Y)$
are based spaces, and there are continuous composition maps
\[ \CM(X,Y) \Smash \CM(Y,Z) \xra{} \CM(X,Z). \]
We again have a strong homotopy category $\hM$, with
$\hM(X,Y)=\pi_0\CM(X,Y)$; when we have defined homotopy groups, we
will also define a weak homotopy category $\hbM$ in the obvious way.

$\CM$ is a closed symmetric monoidal category, with smash product and
function objects again written as $X\Smash Y$ and $F(X,Y)$.  Both of
these constructions are continuous functors of both arguments.  The
unit of the smash product is $S$.

There is a functor $\Sgi\:\CT\xra{}\CM$, such that
\begin{align*}
 \Sgi S^0              &= S                                \\
 \Sgi(A\Smash B)       &= \Sgi A\Smash \Sgi B              \\
 \CM(\Sgi A\Smash X,Y) &= \CT(A,\CM(X,Y))                  \\
 \CM(\Sgi A,\Sgi B)    &= \CT(A,B).
\end{align*}
(For the last of these, see~\cite{el:scw}.)

The last equation shows that $\Sgi$ is a full and faithful embedding
of $\CT$ in $\CM$, so that all of unstable homotopy theory is embedded
in the strong homotopy category $\hM$.  In particular, $\hM$ is very
far from Boardman's stable homotopy category $\CB$.  However, it turns
out that the weak homotopy category $\hbM$ is equivalent to $\CB$.

The definition of this weak homotopy category involves certain
``cofibrant sphere objects'' which we now discuss.  It will be
convenient for us to give a slightly more flexible construction than
that used in~\cite{ekmm:rma}, so as to elucidate certain questions of
naturality.  Let $\CU$ be a universe.  There is a natural way to make
the Lewis-May spectrum $\Sgi\CL(\CU,\rinf)_+$ into a $\LL$-spectrum,
using the action of $\CL(1)=\CL(\rinf,\rinf)$ on $\CL(\CU,\rinf)$ as
well as on the suspension coordinates.  One way to see this is to
observe that $\Sgi\CL(\CU,\rinf)_+=\CL(\CU,\rinf)\ths S^0$, where the
$S^0$ on the right hand side refers to the sphere spectrum indexed on
the universe $\CU$.  We then define
$S(\CU)=S\Smash_\CL\Sgi\CL(\CU,\rinf)_+$.  This gives a contravariant
functor $S\:\{\text{Universes}\}\xra{}\CM$, and it is not hard to
check that $S(\CU)\Smash S(\CV)=S(\CU\oplus\CV)$.

Moreover, for any finite-dimensional subspace $U<\CU$, there is a
natural subobject $S(\CU,U)\leq S(\CU)$ and a canonical isomorphisms
\begin{align*}
 S(\CU,U)\Smash S(\CU,V) &= S(\CU\oplus\CV,U\oplus V)   \\
 \Sg^U S(\CU,U\oplus V)  &= S(\CU,V).
\end{align*}
This indicates that the objects $S(\CU,U)$ are in some sense stable.
They can be defined as follows: take the Lewis-May spectrum
$\Sgi_US^0$ indexed on $\CU$, and then take the twisted half smash
product with the space $\CL(\CU,\rinf)$ to get a Lewis-May spectrum
indexed on $\rinf$ which is easily seen to be an $\LL$-spectrum in a
natural way.  We then apply $S\Smash_\CL(-)$ to get $S(\CU,U)$.

For any $n>0$ and $d\geq 0$ we write
\begin{align*}
 \CL(n)    &\eqdef \CL((\rinf)^n,\rinf)              \\
 S(n)      &\eqdef S((\rinf)^n) = S(1)^{(n)}         \\
 S^d(1)    &\eqdef \Sg^d S(\rinf)                    \\
 S^{-d}(1) &\eqdef S(\rinf,\real^d)
\end{align*}
We will also allow ourselves to write $S^d(n)$ for
$\Sg^kS((\rinf)^n,V)$ where $V$ is a subspace of $(\rinf)^n$ of
dimension $k-d$ and $k$ and $V$ are clear from the context.

Any object of the form $S^U\Smash S(\CV,V)$ is non-canonically
isomorphic to $S^d(1)$, where $d=\dim(U)-\dim(V)$, but when one is
interested in the naturality or otherwise of various constructions it
is often a good idea to forget this fact.  There are isomorphisms
$S^n(1)\Smash S^m(1)\simeq S^{n+m}(1)$ that become canonical and
coherent in the homotopy category.  The homotopy groups of an object
$X\in\CM$ are defined by
\[ \pi_n(X) \eqdef \hM(S^n(1),X). \]
We say that a map $f\:X\xra{}Y$ is a weak equivalence if it induces an
isomorphism $\pi_*(X)\xra{}\pi_*(Y)$, and we define the weak homotopy
category $\hbM$ by inverting weak equivalences.  We define a cell
object to be an object of $\CM$ that is built from the sphere objects
$S^n(1)$ in the usual sort of way; the category $\hbM$ is then
equivalent to the category of cell objects and homotopy classes of
maps.  

\begin{remark}
 In subsequent sections we will consider various spaces of the form
 $\CM(S(1),X)=\Omi F_\CL(S,X)$.  This is weakly equivalent to $\Omi X$
 but not homeomorphic to it; the functor $\Omi\:\CM_S\xra{}\CT$ is not
 representable and has rather poor behaviour.  For this and many
 related reasons it is preferable to replace $X$ by $F_\CL(S,X)$ and
 thus work with EKMM's ``mirror image'' category
 $\CM^S=\{Y\st F_\CL(S,Y)=Y\}$ rather than the equivalent category
 $\CM_S$.  However, our account of these considerations is still in
 preparation so we have used $\CM_S$ in the present work.
\end{remark}

Now let $R$ be a commutative ring object in $\CM$, in other words an
object equipped with maps $S\xra{\eta}R\xla{\mu}R\Smash R$ making the
relevant diagrams geometrically (rather than homotopically)
commutative.  (The term ``ring'' is something of a misnomer, as there
is no addition until we pass to homotopy.)  We let $\CM_R$ denote the
category of module objects over $R$ in the evident sense.  This is
again a topological model category with a closed symmetric monoidal
structure.  The basic cofibrant objects are the free modules
$S^d(1)\Smash R$ for $d\in\Zh$.  The weak homotopy category $\hbM_R$
obtained by inverting weak equivalences is also known as the derived
category of $R$, and written $\DD=\DD_R$; it is equivalent to the
strong homotopy category of cell $R$-modules.  It is not hard to see
that $\DD$ is a monogenic stable homotopy category in the sense
of~\cite{hopast:ash}; in particular, it is a triangulated category
with a compatible closed symmetric monoidal structure.

\section{Strictly unital products}\label{sec-strict-unit}

In the previous sections we worked in the derived category $\DD$ of
(strict) $R$-modules.  In this section we sharpen the picture slightly
by working with modules with strict units.  These are not cell
$R$-modules, so we need to distinguish between
$\hbM_R(X,Y)\eqdef\DD(X,Y)=[X,Y]$ and
$\hM_R(X,Y)\eqdef\pi_0\CM_R(X,Y)=\CM_R(X,Y)/\text{homotopy}$.  Note
that the latter need not have a group structure (let alone an Abelian
one).  However, most of the usual tools of unstable homotopy theory
are available in $\hM_R$, because $\CM_R$ is a topological category
enriched over pointed spaces.  In particular, we will need to use
Puppe sequences.

As previously, we let $x$ be a regular element in $\pi_d(R)$, so $d$
is even.  We regard $x$ as an $R$-module map $S^d(1)\Smash R\xra{}R$,
and we write $R/x$ for the cofibre.  There is thus a pushout diagram
\[ \diagram
  S^d(1)\Smash R \dto_x \rto &
  I\Smash S^d(1)\Smash R \dto \\
  R \rto_\rho & R/x.
\enddiagram 
\]
As $R$ is not a cell $R$-module, the same is true of $R/x$.  However,
the map $\rho\:R\xra{}R/x$ is a $q$-cofibration.  One can also see that
$S^0(1)\Smash R/x$ is a cell $R$-module which is the cofibre in $\DD$
of the map $x\:\Sg^dR\xra{}R$, so it has the homotopy type referred to
as $R/x$ in the previous section.  Moreover, the map 
$S^0(1)\Smash R/x\xra{}R/x$ is a weak equivalence.  It follows that
our new $R/x$ has the same weak homotopy type as in previous sections.

Let $W$ be defined by the following pushout diagram:
\[ \diagram
  R \dto_\rho \rto^\rho &
  R/x \dto^{i_0} \\
  R/x \rto_{i_1} & W
\enddiagram 
\]
There is a unique map $\nabla\:W\xra{}R/x$ such that
$\nabla i_0=1=\nabla i_1$, and there is an evident cofibration
\[ S^{2d+1}(2)\Smash R\xra{}W\xra{}\Rxx. \]
Here 
\[ S^{2d+1}(2) = \Sg S^d(1)\Smash S^d(1) = 
   \begin{cases}
    \Sg^{2d+1}S(\rinf\oplus\rinf)  & \text{ if } d\ge 0 \\
    \Sg S(\rinf\oplus\rinf,\real^{|d|}\oplus\real^{|d|})
      & \text{ if } d<0.
    \end{cases}
\]
We define a \emph{strictly unital product} on $R/x$ to be a map
$\phi\:\Rxx\xra{}R/x$ of $R$-modules such that $\phi|_W=\nabla$.  Let
$P$ be the space of strictly unital products, and let $\ov{P}$ be the
set of products on $R/x$ in the sense of section~\ref{sec-prod-Rx}.

\begin{proposition}\label{prop-strictly-unital}
 The evident map $\pi_0(P)\xra{}\ov{P}$ is a bijection.
\end{proposition}
\begin{proof}
 The cofibration 
 \[ S^{2d+1}(2)\Smash R\xra{i}W\xra{j}\Rxx \]
 gives a fibration $\CM_R(\Rxx,R/x)\xra{j^*}\CM_R(W,R/x)$ of spaces.
 The usual theory of Puppe sequences and fibrations tells us that the
 image of $j^*$ is the union of those components in
 $\pi_0\CM_R(\Rxx,R/x)$ that map to zero in
 $\pi_0\CM_R(S^{2d+1}(2)\Smash R,R/x)=\pi_{2d+1}(R/x)=0$, so $j^*$ is
 surjective.  In particular, we find that $P=(j^*)^{-1}\{\nabla\}$ is
 nonempty.  Similar considerations then show that the $H$-space
 $H=\CM_R(S^{2d+2}(1)\Smash R,R/x)$ acts on $P$, and that for any
 $\phi\in P$ the action map $h\mapsto h.\phi$ gives a weak equivalence
 $H\simeq P$.  This shows that $\pi_0(H)=\pi_{2d+2}(R/x)$ acts freely
 and transitively on $\pi_0(P)$.  This is easily seen to be compatible
 with our free and transitive action of $\pi_{2d+2}(R/x)$ on $\ov{P}$
 (Lemma~\ref{lem-uni-obs}), and the claim follows.
\end{proof}

\begin{remark}\label{rem-associativity}
 These ideas also give another proof of associativity.  Let $Y$ be the
 union of all cells except the top one in $\Rxxx$, so there is a
 cofibration $S^{3d+2}(3)\Smash R\xra{}Y\xra{}\Rxxx$.  Let $\phi$ be a
 product on $R/x$; by the proposition, we may assume that it is
 strictly unital.  It is easy to see that $\phi\circ(\phi\Smash 1)$
 and $\phi\circ(1\Smash\phi)$ have the same restriction to $Y$ (on the
 nose).  It follows using the Puppe sequence that they only differ (up
 to homotopy) by the action of the group $\pi_{3d+3}(R/x)=0$.  Thus,
 $\phi$ is automatically associative up to homotopy.
\end{remark}

We end this section with a more explicit description of the element
$\cb(x)\in\pi_{2d+d}(R)/(2,x)$.  Define $X\eqdef\CM(S^{2d}(2),R/x)$;
this is a space with $\pi_kX=\pi_{2d+k}(R)/x$.  The twist map $\tau$
of $S^{2d}(2)=S^d(1)\Smash S^d(1)$ gives a self-map of $X$, which we
also call $\tau$.  Let $y$ be the map
\[ S^d(1)\Smash S^d(1)\xra{x\Smash x}R\Smash R\xra{\text{mult}}R
   \xra{\rho} R/x,
\]
considered as a point of $X$.  As $R$ is commutative, this is fixed by
$\tau$.  

Next, let $\gm\:I\Smash S^d(1)\xra{}R/x$ be the obvious nullhomotopy
of $x$, and consider the map
\[ I\Smash S^d(1)\Smash S^d(1)\xra{\gm\Smash x}R/x\Smash R
    \xra{\text{mult}} R/x.
\]
This is adjoint to a path $\dl\:I\xra{}X$ with $\dl(0)=0$ and
$\dl(1)=y$.  We could do a similar thing using $x\Smash\gm$ to get
another map $\dl'\:I\xra{}X$, but it is easy to see that
$\dl'=\tau\circ\dl$.  We now define a map
$\phi_0\:\partial(I^2)\xra{}X$ by
\begin{align*}
 \phi_0(s,0) &= 0           \\
 \phi_0(0,t) &= 0           \\
 \phi_0(s,1) &= \dl(s)      \\
 \phi_0(1,t) &= \tau\dl(t).
\end{align*}
We can use the pushout description of $R/x$ to get a pushout
description of $\Rxx$.  Using this, we find that strictly unital
products are just the same as maps $\phi\:I^2\xra{}X$ that extend
$\phi_0$.  

Let $\phi$ be such an extension.  Let $\chi\:I^2\xra{}I^2$ be the
twist map; we find that $\phi'\eqdef\tau\circ\phi\circ\chi$ also
extends $\phi_0$ and corresponds to the opposite product on $R/x$.
Let $U$ be the space $(\{\pm 1\}\tm I^2)/\sim$, where
$(1,s,t)\sim(-1,s,t)$ if $(s,t)\in\partial(I^2)$; clearly this is
homeomorphic to $S^2$.  Define $\psi\:U\xra{}X$ by
$\psi(1,s,t)=\phi(s,t)$ and $\psi(-1,s,t)=\phi'(s,t)=\tau\phi(t,s)$.
It is not hard to see that the class in $\pi_2(X)=\pi_{2d+2}(R)/x$
corresponding to $\psi$ is just $c(\phi)$, and thus that the image in
$\pi_{2d+2}(R)/(2,x)$ is $\cb(x)$.

Another way to think about this is to define a map $\tau\:U\xra{}U$ by
$\tau(r,s,t)=(-r,t,s)$, and to think of $I^2$ as the image of 
$1\tm I^2$ in $U$.  We can then say that $\psi$ is the unique
$\tau$-equivariant extension of $\phi$.

\section{Power operations}\label{sec-pow-op}

In this section, we identify the commutativity obstruction $\cb(x)$ of
Proposition~\ref{prop-Rx} with a kind of power operation.  This is
parallel to a result of Mironov in Baas-Sullivan theory, although the
proofs are independent.  We assume for simplicity that $d\eqdef|x|\geq
0$.

\subsection{The definition of the power operation}\label{subsec-defn-P}

Because $R^*$ is concentrated in even degrees, we know that the
Atiyah-Hirzebruch spectral sequence converging to $R^*\cpi$ collapses
and thus that $R$ is complex orientable.  We choose a complex
orientation once and for all, taking the obvious one if $R$ is (a
localisation of) $MU$.  This gives Thom classes for all complex
bundles.

We write $R^\even(X)$ for the even-degree part of $R^*(X)$, so that
$R^\even(\rp^2)=R^*[\ep]/(2\ep,\ep^2)$.  (In the interesting
applications the ring $R^*$ has no $2$-torsion and so $R^*\rp^2$ has
no odd-degree part.)

We will need notation for various twist maps.  We write $\om$ for the
twist map of $\real^{2d}=\real^d\tm\real^d$, or for anything derived
from that by an obvious functor.  Similarly, we write $s$ for the
twist map of $(\rinf)^2$, and $\sg=S(s)$ for that of 
$S(2)=S(1)\Smash S(1)=S((\rinf)^2)$.  We can thus factor the twist map
$\tau$ of $S^{2d}(2)$ as $\tau=\om\sg=\sg\om$.

We will need to consider the bundle $V^d=\real^{2d}\tm_{C_2}S^2$ over
$S^2/C_2=\rp^2$.  Here $C_2$ is acting on $\real^{2d}$ by $\om$, and
antipodally on $S^2$; the Thom space is $S^2_+\Smash_{C_2}S^{2d}$.  As
$d$ is even, we can regard $V^d$ as $\cplx\ot_\real V^{d/2}$, so we
have a Thom class in $\widetilde{R}^{2d}(S^2_+\Smash_{C_2}S^{2d})$
which generates $\widetilde{R}^\even(S^2_+\Smash_{C_2} S^{2d})$ as a
free module over $R^\even\rp^2=R^*[\ep]/(2\ep,\ep^2)$.

Suppose that $x\in\pi_d(R)$.  Recall that $x$ is represented by a map
$x\:S^d(1)=S^d\Smash S(1)\xra{}R$.  By smashing this with itself and
using the product structure of $R$ we obtain a map
$y\:S^{2d}(2)\xra{}R$.  As $R$ is commutative we have $y\tau=y$.  

Because $S(\CU)$ is a continuous contravariant functor of $\CU$, we
have a map $\CL(2)_+\Smash S(1)\xra{\tht}S(2)$ and thus a map
$\CL(2)_+\Smash S^{2d}\Smash S(1)\xra{y\tht}R$.  If we let
$s\:(\rinf)^2\xra{}(\rinf)^2$ be the twist map and let $C_2$ act on
$\CL(2)$ by $g\mapsto g\circ s$ then $\CL(2)$ is a model for $EC_2$
and thus $\CL(2)/C_2\simeq\rpi$.  As $y\tau=y$ we see that our map
factors through 
$(\CL(2)_+\Smash_{C_2}S^{2d})\Smash S(1)
 \simeq(\rpi)^{V^d}\Smash S(1)$.
For any CW complex $A$, the spectrum $A\Smash S(1)$ is a cofibrant
approximation to $\Sgi A$, so we can regard this map as an element of
$R^0(\rpi)^{V^d}$.  By restricting to $\rp^2$ and using the Thom
isomorphism, we get an element of $R^{-2d}\rp^2$; we define $P(x)$ to
be this element.  We also recall that
$R^\even(\rp^2)=R^*[\ep]/(2\ep,\ep^2)$ and define $\tP(x)$ to be the
coefficient of $\ep$ in $P(x)$, so
$\tP(x)\in{}R^{-2d-2}/2=\pi_{2d+2}(R)/2$.  If $A$ is a CW complex with
only even-dimensional cells then we can replace $R$ by $F(A_+,R)$ to
get power operations $P\:R^{-d}A\xra{}R^{-2d}(\rp^2\tm A)$ and
$\tP\:R^{-d}A\xra{}R^{-2d-2}(A)/2$.  It is not hard to check that this
is the same as the more classical definition given
in~\cite{brmamcst:hir} and thus to deduce the properties listed at the
beginning of Section~\ref{sec-formal}.

We also need a brief remark about the process of restriction to
$\rp^2$.  The space of maps $\mu\:S^2\xra{}\CL(2)$ such that
$\mu(-u)=\mu(u)\circ s$ is easily seen to be contractible.  Choose
such a map $\mu$.  We then have
$(\rp^2)^{V^d}=S^2_+\Smash_{C_2}S^{2d}$, and $P(x)$ is represented by
the composite
\[ (S^2_+\Smash_{C_2}S^{2d})\Smash S(1) \xra{\mu\Smash 1\Smash 1}
   (\CL(2)_+\Smash_{C_2}S^{2d})\Smash S(1) \xra{y\tht} R.
\]
We call this map $\bt_0$.

\subsection{A small modification}\label{subsec-modify}

Let $M$ be the monoid $\CL((\rinf)^2,(\rinf)^2)$.  This acts
contravariantly on $S(2)$, giving a map 
\[ (M_+\Smash_{C_2}S^{2d})\Smash S(2)\xra{} S^{2d}(2) \xra{y} R. \]
Here we use the action of $C_2$ on $M$ given by $g\mapsto g\circ s$.
There is also a homotopically unique map $\lm\:S^2\xra{}M$ such that
$\lm(-u)=\lm(u)\circ s$ for all $u\in S^2$.  By combining this with
the above map, we get a map
\[ (S^2_+\Smash_{C_2}S^{2d})\Smash S(2)\xra{} S^{2d}(2) \xra{y} R. \]
We call this map $\bt_1$.  Recall that in the homotopy category there
is a canonical isomorphism $S(1)\simeq S(2)$, so $\bt_1$ again
represents an element of $R^0(\rp^2)^{V^d}=R^{-2d}\rp^2$.  We claim
that this is the same as $P(x)$.  To see this, choose an isomorphism
$v\:(\rinf)^2\simeq\rinf$, giving a map $v_*\:M\xra{}\CL(2)$ and a map
$S(v)\:S(1)\xra{}S(2)$.  Take $v_*\circ\lm\:S^2\xra{}\CL(2)$ as our
choice of $\mu$, and use $S(v)$ as a representative of the canonical
equivalence $S(1)\simeq S(2)$ in the homotopy category; under these
identifications, $\bt_1$ becomes $\bt_0$.  We leave the rest of the
details to the reader.

\subsection{An adjunction}\label{subsec-adjunction}

Define $Y\eqdef\CM(S^{2d}(2),R)$.  The twist maps $\om$, $\sg$ and
$\tau$ induce commuting involutions of $Y$ with $\tau=\om\sg$.  We can
think of $y\:S^{2d}(2)\xra{}R$ as a point of $Y$, which is fixed under
$\tau$.  The contravariant action of $M$ on $S(2)$ gives a covariant
action on $Y$, which commutes with $\om$.  Using this and our map
$\lm\:S^2\xra{}M$ we can define a map $\bt_2\:S^2_+\xra{}Y$ by
$\bt_2(u)=\lm(u).y$.  If we let $C_2$ act on $Y$ by $\om$ then one
finds that this is equivariant.  We can think of $\bt_2$ as an adjoint
of $\bt_1$ and thus a representative of $P(x)$.

\subsection{Mapping into $R/x$}\label{subsec-Rx}

We are really only interested in the image of $\tP(x)$ in
$\pi_{2d+2}(R)/(2,x)$.  To understand this, we reintroduce the space
$X\eqdef\CM(S^{2d}(2),R/x)$ as in Section~\ref{sec-strict-unit}.  The
unit map $\rho\:R\xra{}R/x$ induces an equivariant map
$\rho_*\:Y\xra{}X$.  We define $\bt_3=\rho_*\circ\bt_2\:S^2_+\xra{}X$.

Note that $\pi_0Y=\pi_{2d}R$, and $\bt_2$ lands in the component
corresponding to $x^2$, and $\rho(x^2)=0$, so $\bt_3$ lands in the
base component.  Moreover, we have $\pi_1X=\pi_{2d+1}(R)/x=0$ and
$C_2$ acts freely on $S^2_+$ so by equivariant obstruction theory we
can extend $\bt_3$ over the cofibre of the inclusion
$S^1_+\xra{}S^2_+$ to get a map $\bt_4$ say.  This cofibre is
equivariantly equivalent to $C_{2+}\Smash S^2$ and 
\[ [C_{2+}\Smash S^2,X]^{C_2}\simeq\pi_2(X)=\pi_{2d+2}(R)/x. \]
It is not hard to see that the element of $\pi_{2d+2}(R)/(2,x)$ coming
from $\bt_4$ is just the image of $\tP(x)$.

\subsection{The abstract argument}\label{subsec-abstract}

We now set up an abstract situation in which we have a space $X$ and
we can define two elements $\al,\bt\in\pi_2(X)/2$ and prove that they
are equal; later we apply this to show that
$\cb(x)=\tP(x)\in\pi_{2d+2}(R)/(2,x)$.  While this involves some
repetition of previous constructions, we believe that it makes the
argument clearer.

Let $M$ be a $2$-connected topological monoid, containing an
involution $\sg$.  Let $C=\{1,\om\}$ be the group of order two, and
define $\tau=\sg\om\in C\tm M$.  Let $X$ be a space with basepoint $0$
and another distinguished point $y$ in the base component.  Suppose
that $C\tm M$ acts on $X$, the whole group fixes $0$, and $\tau$ fixes
$y$.  Suppose also that $\om\simeq 1\:X\xra{}X$ and $\pi_1(X)=0$.

\begin{definition}
 Write $V=\langle\sg,\om\rangle\simeq C_2^2$.  Given $j,k\geq 0$ we
 can let $V$ act on $\real^{j+k}=\real^j\oplus\real^k$ by
 $\om=(-1)\oplus(-1)$ and $\sg=1\oplus(-1)$ so $\tau=(-1)\oplus 1$.
 We write $\real^{j,k}$ for this representation of $V$, and $S^{j,k}$
 for the sphere in $\real^{j,k}$, so that $S^{j,k}=S^{j+k-1}$
 nonequivariantly.
\end{definition}

\begin{definition}\label{defn-alpha}
 Define a $\tau$-equivariant map $\al_0\:S^{0,1}=\{-1,1\}\xra{}X$ by
 $\al_0(-1)=0$ and $\al_0(1)=y$.  Using the evident $\tau$-equivariant
 CW structure on $S^{2,1}$ and the fact that $\pi_1(X)=0$ we find that
 there is an equivariant extension of $\al$ over $S^{2,1}$, which is
 unique modulo $1+\tau^*$.  Nonequivariantly we have $S^{2,1}=S^2$ and
 $\tau\simeq 1\:S^2\xra{}S^2$ so we get a homotopy class of maps
 $S^2\xra{}X$, which is unique modulo $2$.  We write $\al$ for the
 corresponding element of $\pi_2(X)/2$.
\end{definition}

\begin{definition}\label{defn-lambda}
 Define $\lm_0\:S^{0,1}=\{-1,1\}\xra{}M$ by $\lm_0(-1)=1$ and
 $\lm_0(1)=\sg$; this is equivariant with respect to the evident right
 action of $\om$ on $M$.  As $\sg$ acts freely on $S^{2,1}$ and $M$ is
 $2$-connected, we see that there is a $\om$-equivariant extension
 $\lm\:S^{2,1}\xra{}M$, which is unique up to equivariant homotopy.  
\end{definition}

\begin{definition}\label{defn-beta}
 Define $\bt\:S^{2,1}\xra{}X$ by $\bt(u)=\lm(u).y$.  As $\lm$ is
 $\sg$-equivariant and $y$ is fixed by $\tau=\om\sg$ and $\om$
 commutes with $M$ we find that $\bt\sg=\om\bt\:S^{2,1}\xra{}X$.  We
 next claim that $\bt$ can be extended over the cofibre of the
 inclusion of $S^{2,0}$ in $S^{2,1}$ in such a way that we still have
 $\bt\sg=\om\bt$.  This follows easily from the fact that $\sg$ acts
 freely on $S^{2,0}$ and $y$ lies in the base component of $X$ and
 $\pi_1(X)=0$.  The cofibre in question can be identified
 $\sg$-equivariantly with $S^2\Smash\{1,\sg\}_+$.  By composing with
 the inclusion $S^2\xra{}S^2\Smash\{1,\sg\}_+$ we get an element of
 $\pi_2(X)$.  This can be seen to be unique modulo $1+\om_*$ but by
 hypothesis $\om\simeq 1\:X\xra{}X$ so we get a well-defined element
 of $\pi_2(X)/2$, which we also call $\bt$.
\end{definition}

\begin{proposition}\label{prop-alpha-beta}
 We have $\al=\bt\in\pi_2(X)/2$.
\end{proposition}

\begin{proof}
 Consider the following picture of $\real^{2,1}$.
 \[ \epsfbox{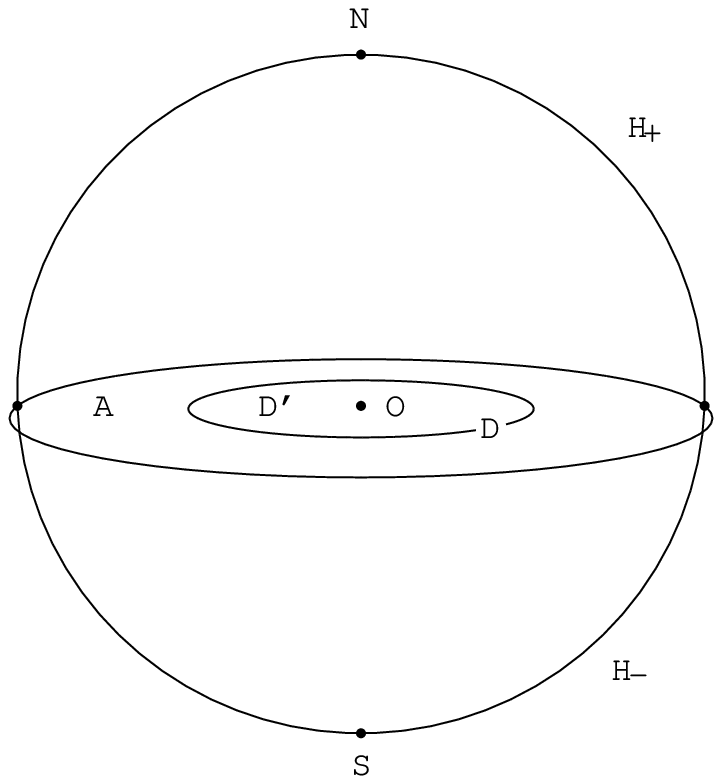} \]
 The axes are set up so that $N=(0,0,1)$ and $S=(0,0,-1)$, so 
 \begin{align*}
  \sg (x,y,z) &= (-x,-y,-z) \\
  \om (x,y,z) &= ( x, y,-z) \\
  \tau(x,y,z) &= (-x,-y, z).
 \end{align*}
 We write $H_+$ and $H_-$ for the upper and lower hemispheres and $D$
 for the unit disc in the plane $z=0$.  Thus $H_+\cup H_-=S^{2,1}$ and
 $H_+\cap H_-=S^{2,0}$, so $H_+\cup H_-\cup D$ can be identified with
 the cofibre of the inclusion $S^{2,0}\xra{}S^{2,1}$.  Note also that
 $H_+\cup D$ is $\tau$-equivariantly homeomorphic to $S^{2,1}$ (by
 radial projection from the $\tau$-fixed point $(0,0,1/2)$, say).

 Let $D'$ be the closed disc of radius $1/2$ centred at $O$ and let
 $A$ be the closure of $D\setminus D'$.  
 Define $\al'\:S^{2,1}\cup D'\xra{}X$ by $\al'=y$ on $S^{2,1}$ and
 $\al'=0$ on $D'$.  We see by obstruction theory that $\al'$ can be
 extended $\tau$-equivariantly over the whole of $S^{2,1}\cup D$.
 Moreover, if we identify $H_+\cup D$ with $S^{2,1}$ as before then
 the restriction of $\al'$ to $H_+\cup D$ represents the same homotopy
 class $\al$ as considered in Definition~\ref{defn-alpha}, as one sees
 directly from the definition.

 Next, note that $S^{2,1}\cup A$ retracts $\sg$-equivariantly onto
 $S^{2,1}$, so we can extend our map $\lm\:S^{2,1}\xra{}M$ over
 $S^{2,1}\cup A$ equivariantly.  As $M$ is $1$-connected, we can
 extend it further over the whole of $S^{2,1}\cup D$, except that we
 have no equivariance on $D'$.

 Now define $\bt'\:S^{2,1}\cup D\xra{}X$ by $\bt'(u)=\lm(u).\al'(u)$.
 We claim that $\bt'\sg=\om\bt'$.  Away from $D'$ this follows easily
 from the equivariance of $\lm$ and $\al'$, and on $D'$ it holds
 because both sides are zero.  Using this and our identification of
 $S^{2,1}\cup D$ with the cofibre of $S^{2,0}\xra{}S^{2,1}$ we see
 that the restriction of $\bt'$ to $H_+\cup D$ represents the class
 $\bt$ in Definition~\ref{defn-beta}.

 Now observe that $S^{2,1}\cup D$ is $2$-dimensional and $M$ is
 $2$-connected, so our map $\lm\:S^{2,1}\cup D\xra{}M$ is
 nonequivariantly homotopic to the constant map with value $1$.  This
 implies that $\al'$ is homotopic to $\bt'$, so $\al=\bt\in\pi_2(X)/2$
 as claimed.
\end{proof}

\subsection{The proof that $\cb(x)=\tP(x)$}
\label{subsec-cb-tP}

We now prove that $\cb(x)=\tP(x)$.  We take
$X\eqdef\CM(S^{2d}(2),R/x)$ and $M\eqdef\CL((\rinf)^2,(\rinf)^2)$ as
before, and define involutions $\om$, $\sg$ and $\tau$ as in
Section~\ref{subsec-defn-P}.  We also define $y$ as in
Section~\ref{subsec-adjunction}.  It is then clear that the map
$\bt_4$ of Section~\ref{subsec-Rx} represents the class $\bt$ of
Definition~\ref{defn-beta}, so that
$\bt=\tP(x)\in\pi_{2d+2}(R)/(2,x)$.  Now consider the constructions at
the end of Section~\ref{sec-strict-unit}.  It is not hard to see that
the space $U$ defined there is $\tau$-equivariantly homeomorphic to
$S^{2,1}$, with the two fixed points being $(0,0,1)$ and $(1,1,1)$.
As the map $\psi\:U\xra{}X$ is equivariant and $\psi(0,0,1)=0$ and
$\psi(1,1,1)=y$, we see that $\psi$ represents the class $\al$ of
Definition~\ref{defn-alpha}, so $\cb(x)=\bt\in\pi_{2d+2}(R)/(2,x)$.
It now follows from Proposition~\ref{prop-alpha-beta} that
$\tP(x)=\cb(x)\pmod{2,x}$, as claimed.


\end{document}